\documentclass[12pt]{article}
\usepackage{amssymb}

\newtheorem{definition}{Definition}[section]
\newtheorem{theorem}[definition]{Theorem}
\newtheorem{lemma}[definition]{Lemma}
\newtheorem{corollary}[definition]{Corollary}
\newtheorem{remark}[definition]{Remark}
\newtheorem{example}[definition]{Example}

\newtheorem{note}[definition]{Note}
\newtheorem{assumption}[definition]{Assumption}
\newtheorem{proposition}[definition]{Proposition}

\def\K{\mathbb K}

\def\K{\mathbb K}

\begin{document}
\title{\bf Two non-nilpotent linear transformations \\
that satisfy
the cubic $q$-Serre relations}
\author{
Tatsuro Ito{\footnote{
Department of Computational Science,
Faculty of Science,
Kanazawa University,
Kakuma-machi,
Kanazawa 920-1192, Japan
}}
$\;$ and
Paul Terwilliger{\footnote{
Department of Mathematics, University of
Wisconsin, 480 Lincoln Drive, Madison WI 53706-1388 USA}
}}
\date{}

\maketitle
\begin{abstract}
Let $\K$ denote an algebraically closed field with
characteristic 0, and 
let $q$ denote a nonzero scalar in $\K$ that is not a root of
unity.
Let ${\mathcal A}_q$ denote the unital associative $\K$-algebra
defined by generators $x,y$ and  relations
\begin{eqnarray*}
x^3y-\lbrack 3\rbrack_q x^2yx
+\lbrack 3\rbrack_q xyx^2
-yx^3&=&0,
\\
y^3x-\lbrack 3\rbrack_q y^2xy
+\lbrack 3\rbrack_q yxy^2
-xy^3&=&0,
\end{eqnarray*}
where
$\lbrack 3 \rbrack_q = (q^3-q^{-3})/(q-q^{-1})$.
We classify up to isomorphism the finite-dimensional
irreducible ${\mathcal A}_q$-modules on which neither 
of $x,y$ is nilpotent.
We discuss how these modules are related to
tridiagonal pairs.

\medskip

 \end{abstract}

\section{Statement of the problem and results}

\noindent
Throughout this paper
 $\K$ will denote an algebraically closed field with
characteristic 0. We fix 
a nonzero scalar
 $q \in \K$ that is not a root of $1$ and
 adopt the
following notation:
\begin{eqnarray}
\lbrack n \rbrack_q = \frac{q^n-q^{-n}}{q-q^{-1}}
\qquad \qquad n = 0,1,2,\ldots 
\label{eq:nbrack}
\end{eqnarray}

\begin{definition}
\label{defa}
\rm
Let ${\mathcal A}_q$ denote the unital associative $\K$-algebra
defined by generators $x,y$
and relations
\begin{eqnarray}
\label{eq:qs1}
x^3y-\lbrack 3\rbrack_q x^2yx
+\lbrack 3\rbrack_q xyx^2
-yx^3&=&0,
\\
\label{eq:qs2}
y^3x-\lbrack 3\rbrack_q y^2xy
+\lbrack 3\rbrack_q yxy^2
-xy^3&=&0.
\end{eqnarray}
We call $x,y$ the {\it standard generators} for ${\mathcal A}_q$.
\end{definition}

\begin{remark}
\rm
The equations 
(\ref{eq:qs1}),
(\ref{eq:qs2})
are the cubic $q$-Serre relations
\cite[p.~11]{lusztig}.
\end{remark}

\noindent
We are interested in a certain class  of 
${\mathcal A}_q$-modules. To describe this class we
recall a concept. Let $V$ denote a 
finite-dimensional vector space over $\K$. 
A linear transformation $X:V\to V$ is said
to be {\it nilpotent} whenever
there exists a positive integer $n$ such that
$X^n=0$. 

\begin{definition}
\label{def:nn}
\rm
Let $V$ denote a finite-dimensional
${\mathcal A}_q$-module.
We say this module is {\it NonNil}
whenever the standard generators $x, y$ are not nilpotent
on $V$.
\end{definition}

\noindent In this paper we will 
classify up to isomorphism the 
NonNil finite-dimensional
irreducible ${\mathcal A}_q$-modules.

\medskip
\noindent Before we state our results, we 
explain why this problem is of interest.
A recent topic of  research
concerns the {\it Leonard
pairs} 
\cite{hartwig},
\cite{nom2},
\cite{Rosengren2},
\cite{LS99},
  \cite{qSerre},
   \cite{LS24},
   \cite{conform},
    \cite{lsint},
\cite{Terint},
\cite{TLT:split},
 \cite{TLT:array},
\cite{qrac},
\cite{aw},
\cite{madrid},
\cite{vidunas}
 and the closely related 
{\it tridiagonal pairs}
\cite{hasan},
\cite{hasan2},
\cite{TD00},
\cite{shape},
\cite{tdanduq},
\cite{nomura},
\cite{nom3},
\cite{nom4}.
A Leonard pair is a pair of semisimple linear transformations
on a finite-dimensional vector space, each of which acts
in an irreducible tridiagonal fashion on an eigenbasis for the other 
\cite[Definition 1.1]{LS99}.
The  Leonard pairs are classified \cite{LS99},
\cite{TLT:array}
and 
correspond
to the orthogonal polynomials that make up the terminating
branch of the Askey scheme 
\cite{KoeSwa},
\cite{qrac}.
A tridiagonal pair
is a mild generalization of a Leonard pair
\cite[Definition 1.1]{TD00}.
For these objects 
the classification problem is open.
Let $V$ denote a NonNil finite-dimensional
irreducible ${\mathcal A}_q$-module. Then the standard generators
$x,y$ act on $V$ as a tridiagonal pair;
see Section 2 for the details.

\medskip
\noindent We now summarize
our classification.
We begin with some comments.
Let $V$ denote a
NonNil
finite-dimensional
irreducible
${\mathcal A}_q$-module. As we will
see, the standard generators $x,y$ are
semisimple on
$V$. Moreover
there exist an integer $d\geq 0$ and nonzero scalars $\alpha, \alpha^* \in \K$
such that the set of distinct eigenvalues of $x$ (resp. $y$) on $V$
is $\lbrace \alpha q^{d-2i} \,|\, 0 \leq i \leq d\rbrace$
(resp. $\lbrace \alpha^*q^{d-2i} \,|\, 0 \leq i \leq d\rbrace$).
We call the ordered pair $(\alpha,\alpha^*)$ the {\it type}
of $V$.
Replacing $x,y$ by $x/\alpha, y/\alpha^*$ the type becomes $(1,1)$.
Consequently it suffices to classify
the NonNil 
finite-dimensional
irreducible
${\mathcal A}_q$-modules of type $(1,1)$.
As we will see, these objects 
are related
to certain modules for the quantum affine algebra
$U_q(\widehat{ \mathfrak{sl}}_2)$.
This algebra is defined as follows.

\begin{definition} 
\label{def:uq}
\rm
\cite[p.~262]{charp}
The quantum affine algebra 
$U_q(\widehat{ \mathfrak{sl}}_2)$ is the unital associative $\K$-algebra 
with
generators $e^{\pm}_i$, $K_i^{{\pm}1}$, $i\in \lbrace 0,1\rbrace $
and the following relations:
\begin{eqnarray}
K_iK^{-1}_i &=& 
K^{-1}_iK_i =  1,
\label{eq:buq1}
\\
K_0K_1&=& K_1K_0,
\label{eq:buq2}
\\
K_ie^{\pm}_iK^{-1}_i &=& q^{{\pm}2}e^{\pm}_i,
\label{eq:buq3}
\\
K_ie^{\pm}_jK^{-1}_i &=& q^{{\mp}2}e^{\pm}_j, \qquad i\not=j,
\label{eq:buq4}
\\
\lbrack e^+_i, e^-_i\rbrack &=& {{K_i-K^{-1}_i}\over {q-q^{-1}}},
\label{eq:buq5}
\\
\lbrack e^{\pm}_0, e^{\mp}_1\rbrack &=& 0,
\label{eq:buq6}
\end{eqnarray}
\begin{eqnarray}
(e^{\pm}_i)^3e^{\pm}_j -  
\lbrack 3 \rbrack_q (e^{\pm}_i)^2e^{\pm}_j e^{\pm}_i 
+\lbrack 3 \rbrack_q e^{\pm}_ie^{\pm}_j (e^{\pm}_i)^2 - 
e^{\pm}_j (e^{\pm}_i)^3 =0, \qquad i\not=j.
\label{eq:buq7}
\end{eqnarray}
We call $e^{\pm}_i$, $K_i^{{\pm}1}$, $i\in \lbrace 0,1\rbrace $
the {\it Chevalley generators} for
$U_q({\widehat{\mathfrak{sl}}}_2)$.
\end{definition}


\begin{remark}
\label{rem:trad}
\rm
\cite[Corollary~3.2.6]{lusztig}
There exists an injection of $\K$-algebras
from ${\mathcal A}_q$ to 
$U_q({\widehat{\mathfrak{sl}}}_2)$ that sends
$x$ and $y$ to $e_0^+$ and 
$e_1^+$ respectively.
Consequently  
 ${\mathcal A}_q$ is often called the
 {\it positive part of
$U_q({\widehat{\mathfrak{sl}}}_2)$}.
\end{remark}

\medskip
\noindent 
The finite-dimensional irreducible
$U_q(\widehat{ \mathfrak{sl}}_2)$-modules
are classified up to isomorphism by
V. Chari and A. Pressley \cite{charp}. In Section
3 we review this classification, and here
 mention just those 
aspects needed to state our results.
Let $\K \lbrack z \rbrack$ denote
the $\K$-algebra consisting
of the polynomials in an indeterminate $z $
that have coefficients in $\K$.
Let $V$ denote a
finite-dimensional 
irreducible
$U_q(\widehat{ \mathfrak{sl}}_2)$-module.
By \cite[Prop.~3.2]{charp}
the actions of $K_0$  and $K_1$ on $V$
are semisimple.
Also by
\cite[Prop.~3.2]{charp}
there exist an integer $d\geq 0$
and scalars $\varepsilon_0, \varepsilon_1$
chosen from $\lbrace 1,-1\rbrace $ such that

\begin{enumerate}
\item the set of distinct eigenvalues of $K_0$ on $V$ is
$\lbrace \varepsilon_0 q^{d-2i} \;|\;0 \leq  i \leq d \rbrace$;
\item $K_0K_1-\varepsilon_0 \varepsilon_1 I$ vanishes
on $V$.
\end{enumerate}
The ordered pair 
$(\varepsilon_0, \varepsilon_1)$ is called
the
{\it type} of $V$.
Now assume $V$ has type $(1,1)$.
Following Drinfel'd
\cite{drinfeld},
Chari and Pressley
define a polynomial $P=P_V$
in 
$\K \lbrack z \rbrack$, called the 
{\it Drinfel'd polynomial of $V$}
\cite[Section~3.4]{charp}.
By \cite[p.~261]{charp}
the map $V\rightarrow P_V$ induces a bijection
between the following two sets:
\begin{enumerate}
\item
 The isomorphism classes of 
finite-dimensional 
irreducible 
$U_q(\widehat{ \mathfrak{sl}}_2)$-modules
of type $(1,1)$;
\item The polynomials  in 
$\K \lbrack z \rbrack$ that have constant coefficient 1.
\end{enumerate}

\medskip
\noindent The main results of this paper
 are contained
in the following two theorems and
subsequent remark.

\begin{theorem}
\label{thm:1}
Let $V$ denote a
NonNil 
finite-dimensional
irreducible
${\mathcal A}_q$-module of type $(1,1)$.
Then there exists
a unique
$U_q(\widehat{ \mathfrak{sl}}_2)$-module
structure on 
$V$ such that the standard generators
$x$ and  $y$ act as $e^+_0+K_0$  
and $e^+_1+K_1$ respectively.
This 
$U_q(\widehat{ \mathfrak{sl}}_2)$-module
is irreducible type $(1,1)$ and $P_V(q^{-1}(q-q^{-1})^{-2})\not=0$.
\end{theorem}

\begin{theorem}
\label{thm:2}
Let $V$ denote a 
finite-dimensional
irreducible 
$U_q(\widehat{ \mathfrak{sl}}_2)$-module of type
$(1,1)$ such that 
$P_V(q^{-1}(q-q^{-1})^{-2})\not=0$.
Then there exists
a unique
${\mathcal A}_q$-module 
structure on 
$V$ such that the standard generators
$x$ and  $y$ act as $e^+_0+K_0$  
and $e^+_1+K_1$ respectively.
This 
${\mathcal A}_q$-module 
is NonNil irreducible 
of type  $(1,1)$.
\end{theorem}

\begin{remark}
\rm Combining 
Theorem
\ref{thm:1}
and Theorem
\ref{thm:2}
we obtain a bijection between the following two sets:
\begin{enumerate}
\item The isomorphism classes of NonNil
finite-dimensional 
irreducible 
${\mathcal A}_q$-modules of type $(1,1)$;
\item The isomorphism classes of 
finite-dimensional 
irreducible 
$U_q(\widehat{ \mathfrak{sl}}_2)$-modules $V$ 
of type $(1,1)$ such that $P_V(q^{-1}(q-q^{-1})^{-2})\not=0$.
\end{enumerate}
\end{remark}

\noindent The plan for the paper is as follows.
In Section 2 we describe how
${\mathcal A}_q$-modules are related to tridiagonal pairs.
In Section 3 we recall some facts about
finite-dimensional irreducible
$U_q(\widehat{ \mathfrak{sl}}_2)$-modules,
including the classification due to Chari and Pressley.
In Section 4 we discuss the Drinfel'd polynomial.
Sections 5--12  are devoted to proving 
Theorem
\ref{thm:1}
and Theorem
\ref{thm:2}.

\section{${\mathcal A}_q$-modules and tridiagonal pairs}

\noindent 
In this section we explain how a NonNil
finite-dimensional irreducible ${\mathcal A}_q$-module
gives a tridiagonal pair.
We will use the following concepts.
Let $V$ denote a finite-dimensional
vector space over $\K$ 
and let $X:V\to V$ denote a linear transformation.
For $\theta \in \K$ we define
\begin{eqnarray*}
V_X(\theta) = \lbrace v \in V \,|\,Xv = \theta v\rbrace.
\end{eqnarray*}
We observe $\theta$ is an eigenvalue of $X$ if and only
if $V_X(\theta)\not=0$, and in this case
$V_X(\theta)$ is the corresponding eigenspace.
Observe the sum
$\sum_{\theta \in \K} V_X(\theta)$ is direct.
Moreover 
this sum is equal to $V$ 
if and only if
$X$ is semisimple.

\begin{lemma}
\label{lem:key}
Let $V$ denote a finite-dimensional
vector space over $\K$.
Let $X:V\to V$ and $Y:V\to V$ denote linear transformations.
Then for all nonzero $\theta \in \K$ the following 
are equivalent:
\begin{enumerate}
\item[{\rm (i)}]  
The expression
$
X^3Y-\lbrack 3\rbrack_q X^2YX
+\lbrack 3\rbrack_q XYX^2
-YX^3
$
vanishes on $V_X(\theta)$.
\item[{\rm (ii)}]
$Y V_X(\theta) \subseteq 
 V_X(q^2\theta)
+
 V_X(\theta)
+
V_X(q^{-2}\theta)
$.
\end{enumerate}
\end{lemma}
\noindent {\it Proof:} 
For $v \in V_X(\theta)$ we
have 
\begin{eqnarray*}
&&(X^3Y-\lbrack 3\rbrack_q X^2YX
+\lbrack 3\rbrack_q XYX^2
-YX^3)v
\\
&& \qquad = 
(X^3-\theta \lbrack 3\rbrack_q X^2
+\theta^2\lbrack 3\rbrack_q X
-\theta^3I)Yv
\qquad \qquad \mbox{Since}\;Xv=\theta v
\\
&&\qquad =
(X-q^2 \theta I)
(X- \theta I)
(X-q^{-2} \theta I)Yv,
\end{eqnarray*}
where $I:V\to V$ is the identity map.
The scalars $q^{2}\theta, \theta, q^{-2}\theta$
are mutually distinct since $\theta\not=0$ and
since $q$ is not a root of $1$.
The result follows.
\hfill $\Box $  

\medskip
\noindent
We now recall the concept of a {\it tridiagonal pair}.

\begin{definition}  
\rm
\cite[Definition~1.1]{TD00}
\label{def:tdp}
Let $V$ denote a vector space over $\K$ with finite positive dimension.
By a {\it tridiagonal pair}  on $V$,
we mean an ordered pair of linear transformations
$A:V\to V$ and $A^*:V \to V$
that satisfy the following four conditions:
\begin{enumerate}
\item Each of $A,A^*$ is semisimple.
\item There exists an ordering $V_0, V_1,\ldots, V_d$ of the  
eigenspaces of $A$ such that 
\begin{equation}
A^* V_i \subseteq V_{i-1} + V_i+ V_{i+1} \qquad \qquad (0 \leq i \leq d),
\label{eq:t1}
\end{equation}
where $V_{-1} = 0$, $V_{d+1}= 0$.
\item There exists an ordering $V^*_0, V^*_1,\ldots, V^*_\delta$ of
the  
eigenspaces of $A^*$ such that 
\begin{equation}
A V^*_i \subseteq V^*_{i-1} + V^*_i+ V^*_{i+1} \qquad \qquad (0 \leq i \leq \delta),
\label{eq:t2}
\end{equation}
where $V^*_{-1} = 0$, $V^*_{\delta+1}= 0$.
\item There does not exist a subspace $W$ of $V$ such  that $AW\subseteq W$,
$A^*W\subseteq W$, $W\not=0$, $W\not=V$.
\end{enumerate}
\end{definition}

\begin{note} \rm
According to a common 
notational convention, $A^*$ denotes the conjugate transpose
of $A$.
 We are not using this convention.
In a tridiagonal pair $A,A^*$ the linear transformations 
$A$ and $A^*$ are arbitrary subject to (i)--(iv) above.
\end{note}

\begin{note} \rm  In Section 1 we mentioned the
concept of a Leonard pair.
A Leonard pair is the same thing 
as a tridiagonal pair for which the $V_i, V^*_i$ all
have dimension 1 \cite[Lemma~2.2]{qSerre}.
\end{note}

\noindent 
We refer the reader to
\cite{TD00},
\cite{qSerre} 
for the basic theory of tridiagonal pairs.
For connections to distance-regular graphs see
\cite[p.~260]{BanIto},
\cite{Cau},
\cite{go},
\cite[Example~1.4]{TD00},
\cite{Leon},
\cite{pasc},
\cite{TersubI}, 
\cite{tersub3}.
For connections to statistical mechanics
see
\cite{bas1}, 
\cite{bas2},
\cite{bas3},
\cite{bas4},
\cite{DateRoan2},
\cite[Section~34]{madrid}  and the
references therein.
For related topics see
\cite{atak1},
\cite{atak2},
\cite{atakdual},
\cite{GYLZmut},
\cite{GYZTwisted},
\cite{GH4},
\cite{GH7},
\cite{GH1},
\cite{cite37},
\cite{koornonn},
\cite{Marco},
\cite{Zhidd},
\cite{ZheCart},
\cite{Zhidden},
\cite{LPcm}.

\medskip
\noindent 
We recall a few basic facts about tridiagonal pairs.
Let $A,A^*$ denote a tridiagonal pair on $V$
and let $d, \delta$ be as in
Definition \ref{def:tdp}(ii), (iii).
 By 
 \cite[Lemma~4.5]{TD00}
 we have
$d=\delta$; we call this common value the {\it diameter}
of $A,A^*$.  
An ordering of the eigenspaces of $A$ (resp. $A^*$)
will be called
{\it standard} whenever it satisfies
(\ref{eq:t1}) (resp. 
(\ref{eq:t2})). 
We comment on the uniqueness of the standard ordering.
Let 
 $V_0, V_1,\ldots, V_d$ denote a standard ordering
of the eigenspaces of $A$.
Then the ordering 
 $V_d, V_{d-1},\ldots, V_0$ 
is standard and no other ordering is standard.
A similar result holds for the 
 eigenspaces of $A^*$.
An ordering of the eigenvalues of $A$ (resp. $A^*$)
will be called {\it standard}
whenever the corresponding ordering of the eigenspaces of $A$
(resp. $A^*$) is standard.

\begin{definition}
\rm
Let $d$ denote a nonnegative integer and let
$\theta_0, \theta_1, \ldots, \theta_d$
denote a sequence of scalars taken from $\K$.
We call this sequence a {\it $q$-string}
whenever there exists a nonzero scalar
$\alpha \in \K$ such that
$\theta_i = \alpha q^{d-2i}$ for $0 \leq i \leq d$. 
\end{definition}

\begin{definition}
\label{def:geom}
\rm
A tridiagonal pair $A,A^*$
is said to be {\it $q$-geometric}
whenever 
(i) there exists a standard ordering of the eigenvalues of
$A$ which forms a $q$-string;
and 
(ii) there exists a standard ordering of the eigenvalues of
$A^*$ which forms a $q$-string.
\end{definition}

\noindent 
We refer the reader to
\cite{hasan},
\cite{hasan2},
\cite{shape},
\cite{tdanduq}, 
\cite{qSerre}
for information about  
$q$-geometric tridiagonal pairs.

\begin{theorem}
{\rm \cite[Lemma~4.8]{qSerre}}  
\label{thm:maintd} 
 Let $V$  denote
a vector space over $\K$ with finite positive dimension.
Let $A:V\to V$ and $A^*:V\to V$
denote linear transformations. Then
 the following are equivalent:
\begin{enumerate}
\item[{\rm (i)}] $A,A^*$ is a $q$-geometric
tridiagonal pair on $V$.
\item[{\rm (ii)}] There exists a
NonNil 
irreducible
${\mathcal A}_q$-module structure on $V$
such that the standard generators $x,y$ act as
$A,A^*$ respectively.
\end{enumerate}
\end{theorem}
\noindent {\it Proof:}
\noindent (i) $\Rightarrow$ (ii):
We first show
\begin{eqnarray}
\label{eq:qs1A}
A^3A^*-\lbrack 3\rbrack_q A^2A^*A
+\lbrack 3\rbrack_q AA^*A^2
-A^*A^3&=&0.
\end{eqnarray}
Since the tridiagonal pair 
$A,A^*$ is
$q$-geometric,
there exists a standard ordering
$\theta_0, \theta_1, \ldots, \theta_d$
of the eigenvalues for $A$ which forms a $q$-string.
For $0 \leq i \leq d$ let
$V_i$ denote the eigenspace of $A$ associated with
$\theta_i$.
Then the ordering $V_0, V_1, \ldots, V_d$ is
 standard and therefore satisfies
(\ref{eq:t1}). By this and since
$\theta_0, \theta_1, \ldots, \theta_d$
is a $q$-string
we find
\begin{eqnarray*}
A^*V_A(\theta) \subseteq 
V_A(q^2\theta)+V_A(\theta)+V_A(q^{-2}\theta)
\end{eqnarray*}
for each eigenvalue $\theta$ of $A$.
Invoking Lemma
\ref{lem:key} we find
that in equation
(\ref{eq:qs1A}) the expression on the left vanishes
on each eigenspace of $A$.
These eigenspaces span $V$ since
 $A$ is semisimple, and equation
(\ref{eq:qs1A}) follows.
Interchanging $A,A^*$ in the above argument
we find
\begin{eqnarray}
\label{eq:qs2A}
A^{*3}A-\lbrack 3\rbrack_q A^{*2}AA^*
+\lbrack 3\rbrack_q A^*AA^{*2}
-AA^{*3}&=&0.
\end{eqnarray}
By 
(\ref{eq:qs1A}), 
(\ref{eq:qs2A}) there exists an ${\mathcal A}_q$-module
structure on $V$ such that 
$x,y$ act as $A,A^*$ respectively.
This
 ${\mathcal A}_q$-module
is irreducible
in view of Definition
\ref{def:tdp}(iv).
This
 ${\mathcal A}_q$-module
is NonNil since 
each of $A,A^*$ has all eigenvalues  nonzero.
\\
\noindent (ii) $\Rightarrow $ (i):
The field $\K$ is algebraically closed so it
contains all the eigenvalues
for $A$.
Since $V$ has finite positive dimension, 
$A$ has at least one eigenvalue.
Since 
$x$ acts on $V$ as $A$ and since
this action is 
not nilpotent,
 $A$
has at 
least one nonzero 
eigenvalue $\theta$.
The scalars $\theta, \theta q^{2}, \theta q^{4},\ldots$
are mutually distinct since $q$ is not a root of unity,
so these scalars are not all
eigenvalues for $A$.
Consequently there exists a nonzero
$\eta \in \K$
such that $V_A(\eta)\not=0$ and
$V_A(\eta q^{2})=0$.
There exists an integer $d\geq 0$ such
that
$V_A(\eta q^{-2i})$  is nonzero
for $0 \leq i \leq d$ and zero for $i=d+1$.
We abbreviate $V_i = V_A(\eta q^{-2i})$ for $0 \leq i \leq d$. 
By  construction $\sum_{i=0}^d V_i$ is 
$A$-invariant.
By 
Lemma
\ref{lem:key} we find
\begin{eqnarray*}
A^*V_i \subseteq V_{i-1}+V_i + V_{i+1}
\qquad \qquad (0 \leq i \leq d),
\end{eqnarray*}
where $V_{-1}=0$, $V_{d+1}=0$.
Consequently 
$\sum_{i=0}^d V_i$ is 
$A^*$-invariant.
Since the ${\mathcal A}_q$-module
 structure
on $V$ is irreducible and since
$\sum_{i=0}^d V_i\not=0$
we find
$\sum_{i=0}^d V_i=V$.
Now $A$ is semisimple and
Definition \ref{def:tdp}(ii) holds.
Interchanging $A,A^*$ in the above 
argument we find
$A^*$ is semisimple  and
Definition \ref{def:tdp}(iii) holds.
Definition \ref{def:tdp}(iv) is satisfied
since the ${\mathcal A}_q$-module
$V$ is irreducible.
We have now shown that
$A,A^*$ is a tridiagonal pair on $V$.
From the construction
this tridiagonal pair is $q$-geometric.
\hfill $\Box $


\begin{corollary}
\label{cor:maintd} 
Let $V$ denote a  NonNil
finite-dimensional irreducible
${\mathcal A}_q$-module.
Then the standard generators $x,y$ are semisimple on $V$.
 Moreover
there exist an integer $d\geq 0$ and nonzero scalars $\alpha, \alpha^* \in \K$
such that the set of distinct eigenvalues of $x$ (resp. $y$) on $V$
is $\lbrace \alpha q^{d-2i} \,|\, 0 \leq i \leq d\rbrace$
(resp. $\lbrace \alpha^*q^{d-2i} \,|\, 0 \leq i \leq d\rbrace$).
\end{corollary}
\noindent {\it Proof:}  
Immediate from
Definition
\ref{def:geom} and
Theorem
\ref{thm:maintd}.
\hfill $\Box $  

\begin{definition}
\label{def:nntype}
\rm
Let $V$ denote a  NonNil
finite-dimensional irreducible
${\mathcal A}_q$-module.
By the {\it type} of $V$ we mean
the ordered pair $(\alpha, \alpha^*)$ from
Corollary
\ref{cor:maintd}.
By the {\it diameter} of $V$ we mean the scalar
$d$ from
Corollary
\ref{cor:maintd}.
\end{definition}

\begin{lemma}
\label{lem:autuq2}
For all nonzero $\alpha, \alpha^* \in \K$
there exists a $\K$-algebra
automorphism of 
${\mathcal A}_q$ such that
\begin{eqnarray*}
x \rightarrow \alpha x,
\qquad \qquad y\rightarrow \alpha^*y.
\end{eqnarray*}
\end{lemma}
\noindent {\it Proof:}  This is immediate from
Definition
\ref{defa}.
\hfill $\Box $  

\begin{remark}
\rm   
Given a NonNil 
 finite-dimensional irreducible
${\mathcal A}_q$-module,
we can change its 
 type to any other type
by  applying an automorphism from 
Lemma
\ref{lem:autuq2}.
\end{remark}

\section{Finite dimensional
$U_q(\widehat{ \mathfrak{sl}}_2)$-modules}

\noindent
In this section we recall
some facts about 
finite-dimensional
irreducible $U_q(\widehat{ \mathfrak{sl}}_2)$-modules,
including the classification due to
Chari and Pressley \cite{charp}.

\medskip
\noindent We begin with some notation. 
Let $V$ denote a vector space over $\K$ with finite
positive dimension.
Let $d$ denote a nonnegative integer.
By a {\it decomposition of $V$
with diameter $d$} we mean a sequence
$U_0, U_1, \ldots, U_d$
consisting 
of nonzero subspaces of $V$ such that
\begin{eqnarray*}
V = \sum_{i=0}^d U_i \qquad \qquad (\mbox{direct sum}).
\end{eqnarray*}
We do not assume the spaces $U_0,U_1, \ldots, U_d$ have dimension 1.
For notational convenience we define $U_{-1}=0$ and $U_{d+1}=0$.

\begin{lemma}\label{lem:nnweight}
{\rm \cite[Prop.~3.2]{charp}}  
Let $V$ denote a finite-dimensional irreducible 
$U_q(\widehat{ \mathfrak{sl}}_2)$-module.
Then there exist
unique scalars
$\varepsilon_0,
\varepsilon_1$ in 
$\lbrace 1,-1\rbrace $
and   a unique decomposition
$U_0, \ldots, U_d$ of $V$ 
such that  
\begin{eqnarray}\label{eq:nnkmove}
(K_0-\varepsilon_0q^{2i-d}I)U_i=0,
\qquad \qquad  
 (K_1-\varepsilon_1q^{d-2i}I)U_i=0
\end{eqnarray}
for all $i=0,1,\ldots, d$.
Moreover,  for $0 \leq i \leq d$ we have 
\begin{eqnarray}
&&
e_0^+
U_i \subseteq U_{i+1}, \qquad
e_1^-
U_i \subseteq U_{i+1}, 
\label{eq:altemove1}
\\
&&
e_0^-
U_i \subseteq U_{i-1}, \qquad
e_1^+
U_i \subseteq U_{i-1}.
\label{eq:altemove2}
\end{eqnarray}
\end{lemma}

\begin{definition}
\label{def:type}
\rm   
 The ordered pair 
 $(\varepsilon_0, \varepsilon_1)$ in Lemma  
\ref{lem:nnweight} is
  the {\it type}
 of $V$ and   $d$ is the {\it diameter} of $V$.
The sequence 
$U_0,\ldots,U_d$
is the 
{\it weight space} decomposition
of
 $V$ (relative to $K_0$ and $K_1$). 
\end{definition} 

\begin{lemma}
\label{lem:autuq}
{\rm \cite[Prop.~3.3]{charp}}
For any choice of scalars $\varepsilon_0, \varepsilon_1$ 
{f}rom $\lbrace 1,-1\rbrace $,
there exists a $\K$-algebra
automorphism of 
$U_q(\widehat{ \mathfrak{sl}}_2)$
such that
\begin{eqnarray*}
K_i \rightarrow \varepsilon_i K_i,
\qquad \qquad e^+_i\rightarrow \varepsilon_i e^+_i,
\qquad \qquad e^-_i\rightarrow  e^-_i
\end{eqnarray*}
for $i \in \lbrace 0,1\rbrace $.
\end{lemma}

\begin{remark} 
\label{rem:alter}
\rm
Given a finite-dimensional irreducible
$U_q(\widehat{ \mathfrak{sl}}_2)$-module,
we can alter its type to any other type
by  applying an automorphism from Lemma  \ref{lem:autuq}.
\end{remark}

\noindent We now recall the evaluation modules.

\begin{lemma}
\label{lem:evmod}
{\rm \cite[Section~4]{charp}}
There exists a family of finite-dimensional
irreducible type $(1,1)$
$U_q(\widehat{ \mathfrak{sl}}_2)$-modules 
\begin{eqnarray*}
V(d,a) \qquad \qquad d=1,2,\ldots \qquad \qquad 0 \neq a \in \K
\end{eqnarray*}
with the following property.
The module $V(d,a)$ has a basis $v_0, v_1, \ldots, v_d$ such
that
\begin{eqnarray*}
K_0 v_i &=& q^{2i-d}v_i \qquad \qquad (0 \leq i \leq d),
\\
e_0^- v_i &=& qa^{-1}\lbrack d-i+1 \rbrack_qv_{i-1}\qquad (1 \leq i \leq d),
\qquad e_0^- v_0 = 0,
\\
e_0^+ v_i &=& q^{-1}a\lbrack i+1 \rbrack_qv_{i+1}\qquad (0 \leq i \leq d-1),
\qquad e_0^+ v_d = 0,
\\
K_1 v_i &=& q^{d-2i}v_i \qquad \qquad (0 \leq i \leq d),
\\
e_1^- v_i &=& \lbrack i+1 \rbrack_qv_{i+1}\qquad (0 \leq i \leq d-1),
\qquad e_1^- v_d = 0,
\\
e_1^+ v_i &=& \lbrack d-i+1 \rbrack_qv_{i-1}\qquad (1 \leq i \leq d),
\qquad e_1^+ v_0 = 0.
\end{eqnarray*}
\end{lemma}

\begin{definition}
\rm
\label{def:ev}
\cite[Definition~4.2]{charp}
The 
$U_q(\widehat{ \mathfrak{sl}}_2)$-module
$V(d,a)$ from Lemma
\ref{lem:evmod} is called an {\it evaluation module}.
The scalar $a$ is called the {\it evaluation parameter}.
\end{definition}

\begin{remark}
\rm
\cite[Remark~4.2]{charp}
Referring to
Lemma \ref{lem:evmod} and
Definition \ref{def:ev},
two evaluation modules
$V(d,a)$ and $V(d',a')$ are isomorphic
if and only if $(d,a)=(d',a')$.
\end{remark}

\noindent 
We now recall how the tensor product of
two 
$U_q(\widehat{ \mathfrak{sl}}_2)$-modules becomes
 a 
$U_q(\widehat{ \mathfrak{sl}}_2)$-module.
We will use the following fact.

\begin{lemma}
{\rm \cite[p.~263]{charp}}
\label{lem:hopf} 
$U_q({\widehat{\mathfrak{sl}}}_2)$ has the following Hopf algebra
structure. The comultiplication $\Delta$ satisfies
\begin{eqnarray*}
\Delta(e_i^{+})&=& e_i^+\otimes K_i+1 \otimes e_i^+,
\\
\Delta(e_i^{-})&=& e_i^-\otimes 1+ K_i^{-1}\otimes e_i^-,
\\
\Delta(K_i)&=& K_i\otimes K_i.
\end{eqnarray*}
The counit $\varepsilon$ satisfies
\begin{eqnarray*}
\varepsilon(e_i^{\pm})=0, \qquad \qquad
\varepsilon(K_i)=1.
\end{eqnarray*}
The antipode $S$ satifies
\begin{eqnarray*}
S(K_i)= K_i^{-1}, 
\qquad 
\qquad
S(e_i^+)= -e_i^+K_i^{-1}, 
\qquad 
\qquad
S(e_i^-)= -K_ie_i^-.
\end{eqnarray*}
\end{lemma}

\noindent Combining Lemma \ref{lem:hopf}
with 
\cite[p.~110]{charpbook} we routinely obtain the following.

\begin{lemma}
\label{lem:hopfmod}
Let $V,W$ denote 
$U_q(\widehat{ \mathfrak{sl}}_2)$-modules.
Then the tensor product
$V\otimes W=V\otimes_\K W$ has  
the following 
$U_q(\widehat{ \mathfrak{sl}}_2)$-module 
structure. 
For $v \in V$, for $w \in W$ and for $i \in \lbrace 0,1\rbrace $, 
\begin{eqnarray*}
e_i^{+}(v\otimes w)&=& e_i^+v\otimes K_iw+v \otimes e_i^+w,
\\
e_i^{-}(v\otimes w)&=& e_i^-v\otimes w+ K_i^{-1}v\otimes e_i^-w,
\\
K_i(v\otimes w)&=& K_iv\otimes K_iw.
\end{eqnarray*}
\end{lemma}

\begin{definition}
\rm
\cite[p.~110]{charpbook}
There exists
a one dimensional
$U_q(\widehat{ \mathfrak{sl}}_2)$-module on which
each element $z \in
U_q(\widehat{ \mathfrak{sl}}_2)$
acts as $\varepsilon(z)I$,
where $\varepsilon$ is
from 
Lemma 
\ref{lem:hopf} and $I$ is the identity map.
In particular on this module 
each of
$e^{\pm}_0$,
$e^{\pm}_1$
vanishes
and each of
$K^{\pm 1}_0$, 
$K^{\pm 1}_1$ 
acts as $I$. This module
is irreducible, with diameter $0$ and type $(1,1)$.
This module is unique up
to isomorphism.
We call this module
the
 {\it trivial}
$U_q(\widehat{ \mathfrak{sl}}_2)$-module.
\end{definition}

\noindent 
We now state
Chari and Pressley's 
classification of the finite-dimensional irreducible
$U_q(\widehat{ \mathfrak{sl}}_2)$-modules.
In view of
Remark \ref{rem:alter}
it suffices to consider the modules
of type $(1,1)$.

\begin{theorem}
\label{thm:class}
{\rm \cite[Theorem~4.8, Theorem~4.11]{charp}}
Each nontrivial finite-dimensional irreducible
$U_q(\widehat{ \mathfrak{sl}}_2)$-module of type $(1,1)$ 
is isomorphic to a tensor product of evaluation modules.
Two such tensor products are isomorphic if and only if
one is obtained from the other by permuting the factors
in the tensor product.
A tensor product of evaluation modules
\begin{eqnarray*}
V(d_1, a_1)\otimes 
V(d_2, a_2)\otimes 
\cdots 
\otimes
V(d_N, a_N)
\end{eqnarray*}
is irreducible if and only if
\begin{eqnarray*}
a_i a_j^{-1} \not\in \lbrace
q^{d_i+d_j},
q^{d_i+d_j-2},
\ldots, 
q^{|d_i-d_j|+2}\rbrace,
\qquad \qquad (1 \leq i,j\leq N, \; i\not =j)
\end{eqnarray*}
and in this case it is type $(1,1)$.
\end{theorem}

\noindent We end this section with a comment.


\begin{lemma}
\label{cor:genf}
Let $V$ denote a nontrivial finite-dimensional irreducible
$U_q(\widehat{ \mathfrak{sl}}_2)$-module of type $(1,1)$.
Write $V$ as a tensor product of evaluation modules:
\begin{eqnarray*}
V = V(d_1, a_1) \otimes
V(d_2, a_2) \otimes
\cdots \otimes
V(d_N, a_N).
\end{eqnarray*}
Let $U_0, \ldots, U_d$ denote the
weight space decomposition of $V$ from 
Lemma
\ref{lem:nnweight}. Then
\begin{eqnarray*}
\sum_{i=0}^d \mbox{dim}(U_i)z^i=
\prod_{j=1}^N
(1+z+z^2+\cdots + z^{d_j}). 
\end{eqnarray*}
\end{lemma}
\noindent {\it Proof:} 
Combine
Lemma \ref{lem:evmod}
and 
Lemma
\ref{lem:hopfmod}. 
\hfill $\Box $ \\

\section{The Drinfel'd polynomial}

\noindent In this section we recall the Drinfel'd polynomial
 associated with a finite-dimensional
 irreducible 
$U_q(\widehat{ \mathfrak{sl}}_2)$-module
of type $(1,1)$.

\begin{definition}
\label{def:sig}
\rm
Let $V$ denote a finite-dimensional irreducible 
$U_q(\widehat{ \mathfrak{sl}}_2)$-module
of type $(1,1)$, and let $U_0, \ldots, U_d$
denote the corresponding weight space decomposition
from Lemma
\ref{lem:nnweight}.
From Lemma  
\ref{cor:genf} we find
 $\mbox{dim}(U_0)=1$.
Pick an integer $i \geq 0$.
By Lemma
\ref{lem:nnweight} we find
$U_0$ is an eigenspace
for $(e^+_1)^i(e^+_0)^i$; let $\sigma_i$ denote
the corresponding eigenvalue.
We observe
$\sigma_i=0$ if $i>d$.
\end{definition}

\noindent We will use the following notation.
With reference to
(\ref{eq:nbrack}) we define
\begin{eqnarray}
\lbrack n \rbrack^!_q = 
\lbrack n \rbrack_q
\lbrack n-1 \rbrack_q
\cdots 
\lbrack 1 \rbrack_q
\qquad \qquad n=0,1,2,\ldots
\label{eq:nfac}
\end{eqnarray}
We interpret 
$\lbrack 0 \rbrack^!_q =1$. 

\begin{definition}
\rm
\cite[Section~3.4]{charp}
\label{def:pv}
Let $V$ denote a finite-dimensional irreducible
$U_q(\widehat{ \mathfrak{sl}}_2)$-module
of type $(1,1)$. We
define a polynomial $P_V
\in 
 \K \lbrack z \rbrack$ by 
\begin{eqnarray}
P_V = \sum_{i=0}^\infty \frac{(-1)^i \sigma_i q^i z^i}{
\lbrack i \rbrack^{!\,2}_q},
\label{eq:pv}
\end{eqnarray}
where the scalars $\sigma_i $ are from
Definition 
\ref{def:sig}.
We observe that $P_V$ has degree at most 
the diameter of $V$. Moreover $P_V$ has
constant coefficent $\sigma_0=1$.
We call $P_V$
the {\it Drinfel'd polynomial of $V$}.
\end{definition}

\begin{note}
\rm
The polynomial $P_V$ from
Definition
\ref{def:pv} is the same as the polynomial
$P$ from
\cite[Theorem~3.4]{charp}, provided
$P$ is normalized to have constant coefficient
1. This is explained in 
\cite[p.~268]{charp}.
\end{note}

\noindent The Drinfel'd polynomial has the following
properties.

\begin{theorem}
{\rm \cite[Corollary~4.2, Proposition~4.3]{charp}}
Let $V$ denote a finite-dimensional irreducible
$U_q(\widehat{ \mathfrak{sl}}_2)$-module of type $(1,1)$.
If $V=V(d,a)$ is an evaluation module
then 
\begin{eqnarray*}
P_V = (1-q^{d-1}az)
 (1-q^{d-3}az)
\cdots
 (1-q^{1-d}az).
\end{eqnarray*}
If $V=V_1 \otimes \cdots \otimes
V_N$ is a tensor product of evaluation modules
then
$P_V = \prod_{i=1}^N P_{V_i}$. 
\end{theorem}

\begin{theorem}
\label{thm:drinbij}
{\rm \cite[p.~261]{charp}}
The map $V\rightarrow P_V$ induces a bijection
between the following two sets:
\begin{enumerate}
\item[{\rm (i)}]  
 The isomorphism classes of 
finite-dimensional 
irreducible 
$U_q(\widehat{ \mathfrak{sl}}_2)$-modules
of type $(1,1)$;
\item[{\rm (ii)}]  
The polynomials  in 
$\K \lbrack z\rbrack$ that have constant coefficient 1.
\end{enumerate}
\end{theorem}

\section{
From 
NonNil ${\mathcal A}_q$-modules to
$U_q(\widehat{ \mathfrak{sl}}_2)$-modules}

\noindent In this section we cite 
some results from 
 \cite{tdanduq} that will be of use
in our proof of
Theorem
\ref{thm:1}.

\medskip
\noindent 
The following theorem is
a minor adaption 
\cite[Theorem~3.3]{tdanduq} and \cite[Theorem~13.1]{tdanduq};
we will sketch the proof in order to motivate
what comes later in the paper.

\begin{theorem}
\label{thm:towardsuq}
{\rm \cite[Theorem~3.3, Theorem~13.1]{tdanduq}}
Let $V$ denote a
NonNil 
finite-dimensional
irreducible
${\mathcal A}_q$-module of type $(1,1)$.
Then there exists
a unique
$U_q(\widehat{ \mathfrak{sl}}_2)$-module
structure on 
$V$ such that the standard generators
$x$ and $y$ act as $e^+_0+K_0$  
and $e^+_1+K_1$ respectively.
This 
$U_q(\widehat{ \mathfrak{sl}}_2)$-module
is irreducible and type $(1,1)$.
\end{theorem}
\noindent {\it Sketch of Proof:} 
For notational convenience let $A:V\to V$
and $A^*:V\to V$ denote the linear transformations
such that $x$ and $y$ act on $V$ as $A$ and $A^*$ respectively. 
By Theorem
\ref{thm:maintd} the pair
 $A,A^*$ is a $q$-geometric tridiagonal pair on $V$.
Let $d$ denote the diameter of this tridiagonal pair.
By Definition
\ref{def:geom} and since the
${\mathcal A}_q$-module
$V$ has type $(1,1)$,
the sequence
$q^{d}, q^{d-2}, \ldots, q^{-d}$
is a standard ordering of the eigenvalues 
for each of $A, A^*$.
For $0 \leq i \leq d$ let
 $V_i$ (resp. $V^*_i$) denote the eigenspace
 of $A$ (resp. $A^*$) associated with the eigenvalue
$q^{2i-d}$ (resp. $q^{d-2i}$).
Moreover 
for $0 \leq i\leq d$ we define
\begin{eqnarray*}
U_i &=&(V^*_0+\cdots + V^*_i)\cap (V_i + \cdots + V_d).
\end{eqnarray*}
By \cite[Lemma~4.2]{tdanduq}
the sequence $U_0,\ldots, U_d$
is a decomposition of $V$.
Therefore there exists a linear transformation $K:V\to V$
such that for $0 \leq i \leq d$, $U_i$ is an eigenspace for
$K$ with eigenvalue $q^{2i-d}$.
For $0 \leq i\leq d$ we define
\begin{eqnarray*}
W_i &=&(V^*_{0}+\cdots + V^*_i)\cap (V_0 + \cdots + V_{d-i}),
\\
W^*_i &=&(V^*_{d-i}+\cdots + V^*_d)\cap (V_i + \cdots + V_d).
\end{eqnarray*}
Then 
$W_0,\ldots, W_d$ and 
$W^*_0,\ldots, W^*_d$ are decompositions of $V$
\cite[Lemma~4.2]{tdanduq}.
Therefore there exist linear transformations
$B:V\to V$ and $B^*:V\to V$ such that
for $0 \leq i \leq d$, $W_i$ (resp. $W^*_i$) is an eigenspace
for $B$ (resp. $B^*$) with eigenvalue $q^{2i-d}$ (resp. $q^{d-2i}$).
By \cite[Theorem~7.1]{tdanduq} we have
\begin{eqnarray}
\label{eq:start}
\frac{qAB-q^{-1}BA}{q-q^{-1}}&=&I,
\\
\frac{qBA^*-q^{-1}A^*B}{q-q^{-1}}&=&I,
\\
\frac{qA^*B^*-q^{-1}B^*A^*}{q-q^{-1}}&=&I,
\\
\frac{qB^*A-q^{-1}AB^*}{q-q^{-1}}&=&I
\end{eqnarray}
and by \cite[Theorem~10.1]{tdanduq} we have
\begin{eqnarray}
\frac{qK^{-1}A-q^{-1}AK^{-1}}{q-q^{-1}}&=&I,
\\
\frac{qBK^{-1}-q^{-1}K^{-1}B}{q-q^{-1}}&=&I,
\\
\frac{qKA^*-q^{-1}A^*K}{q-q^{-1}}&=&I,
\\
\frac{qB^*K-q^{-1}KB^*}{q-q^{-1}}&=&I.
\end{eqnarray}
By construction we have
\begin{eqnarray}
A^3A^*-
\lbrack 3\rbrack_q A^2A^*A
+\lbrack 3\rbrack_q AA^*A^2
-A^*A^3&=&0,
\\
A^{*3}A-\lbrack 3\rbrack_q A^{*2}AA^*
+\lbrack 3\rbrack_q A^*AA^{*2}
-AA^{*3}&=&0
\end{eqnarray}
and by \cite[Theorem~12.1]{tdanduq} we have
\begin{eqnarray}
B^3B^*-
\lbrack 3\rbrack_q B^2B^*B
+\lbrack 3\rbrack_q BB^*B^2
-B^*B^3&=&0,
\\
B^{*3}B-\lbrack 3\rbrack_q B^{*2}BB^*
+\lbrack 3\rbrack_q B^*BB^{*2}
-BB^{*3}&=&0.
\label{eq:fin}
\end{eqnarray}

\noindent In order to connect the above equations with the
defining relations for 
$U_q(\widehat{ \mathfrak{sl}}_2)$ we make a change of variables.
Define 
\begin{eqnarray*}
R&=&A-K,
\qquad \qquad 
\qquad 
L=A^*-K^{-1},\\
r&=& \frac{I-KB^*}{q(q-q^{-1})^2},
\qquad \qquad 
l= \frac{I-K^{-1}B}{q(q-q^{-1})^2},
\label{eq:ldef}
\end{eqnarray*}
so that
\begin{eqnarray}
A&=&K+R,
\qquad \qquad \quad\qquad \qquad
A^*=K^{-1}+L,
\label{eq:ared}
\\
B&=& K-q(q-q^{-1})^2Kl,
\qquad \qquad 
B^*= K^{-1}-q(q-q^{-1})^2K^{-1}r.
\label{eq:bsred}
\end{eqnarray}
Evaluating 
(\ref{eq:start})--(\ref{eq:fin}) using
(\ref{eq:ared}), (\ref{eq:bsred}) we routinely
obtain
\begin{eqnarray*}
&&K R K^{-1} = q^2 R,
\qquad  \qquad 
K L K^{-1}= q^{-2} L,
\label{eq:app1}
\\
&&
K r K^{-1} = q^2 r, \qquad \qquad  
 K l K^{-1} = q^{-2} l,  
\label{eq:appr}
\\
&&
rL-Lr = {{K-K^{-1}}\over {q-q^{-1}}},
\qquad \qquad 
lR-Rl = {{K^{-1}-K}\over {q-q^{-1}}},
\label{eq:app3}
\\
&&
lL = Ll , \qquad \qquad rR=Rr,
\label{eq:app4}
\\
&&
0= R^3L - \lbrack 3 \rbrack_q R^2LR+ \lbrack 3 \rbrack_q RLR^2 - LR^3, 
\label{eq:app5}
\\
&&
0= L^3R - \lbrack 3 \rbrack_q L^2RL+ \lbrack 3 \rbrack_q LRL^2 - RL^3.
\label{eq:app6}
\\
&&
0= r^3l - \lbrack 3 \rbrack_q r^2lr+ \lbrack 3 \rbrack_q rlr^2 - lr^3, 
\label{eq:app7}
\\
&&
0= l^3r - \lbrack 3 \rbrack_q l^2rl+ \lbrack 3 \rbrack_q lrl^2 - rl^3.
\label{eq:app8}
\end{eqnarray*}
Consequently $V$ becomes
a $U_q(\widehat{\mathfrak{sl}}_2)$-module
on which
the Chevalley generators act as follows:

\medskip
\centerline{
\begin{tabular}[t]{c|cccccccc}
        {\rm generator}  
       	& $e_0^+$  
         & $e_1^+$  
         & $e_0^-$  
         & $e_1^-$  
         & $K_0$  
         & $K_1$  
         & $K_0^{-1}$  
         & $K_1^{-1}$  
	\\
	\hline 
{\rm action on $V$} 
&  $R$ & $L$ & $l$ & $r$ & $K$ & $K^{-1}$ & $K^{-1}$ & $K$ 
\end{tabular}}

\medskip
\noindent
By the construction
$x$ and $y$ act on $V$ as $e^+_0+K_0$ and
$e^+_1+K_1$ respectively.
By this and since
the ${\mathcal A}_q$-module $V$
 is irreducible 
we find
the
$U_q(\widehat{\mathfrak{sl}}_2)$-module $V$ is
irreducible.
By the definition of $K$ and since
$K_0, K_1$ act on $V$ as $K,K^{-1}$ respectively
we find
the $U_q(\widehat{\mathfrak{sl}}_2)$-module $V$ is type $(1,1)$.
We have now proved
all the assertions of the theorem except uniqueness.
The proof of uniqueness is given in
\cite[Section~14]{tdanduq}.
\hfill $\Box $ \\

\begin{remark}
\rm
The equations (\ref{eq:start})--(\ref{eq:fin})
give essentially the equitable presentation of
$U_q(\widehat{ \mathfrak{sl}}_2)$
\cite{tdanduq},
\cite{equit1}, \cite{equit2}.
\end{remark}

\section{From 
$U_q(\widehat{ \mathfrak{sl}}_2)$-modules to 
NonNil ${\mathcal A}_q$-modules, I}

\noindent 
Let $V$ denote a finite-dimensional irreducible
$U_q(\widehat{ \mathfrak{sl}}_2)$-module of type $(1,1)$.
In this section we show 
that there exists 
a NonNil  
${\mathcal A}_q$-module structure
on $V$ such that
the standard generators $x$ and $y$ act as
$e^+_0+K_0$ and
$e^+_1+K_1$ respectively.
We discuss some of the basic properties of
this 
${\mathcal A}_q$-module.

\begin{lemma}
\label{lem:uqqs}
There exists a homomorphism of
$\K$-algebras from
${\mathcal A}_q$
to 
$U_q(\widehat{ \mathfrak{sl}}_2)$
that sends $x$ and $y$ to
 $e^+_0+K_0$  
and $e^+_1+K_1$ respectively.
\end{lemma}
\noindent {\it Proof:} 
Using
the defining relations
for 
$U_q(\widehat{ \mathfrak{sl}}_2)$ we routinely
find
 that 
 $e^+_0+K_0$  
and $e^+_1+K_1$ satisfy the 
cubic $q$-Serre relations.
\hfill $\Box $ \\

\begin{assumption}
\label{as}
\rm Throughout this section $V$ will
denote a finite-dimensional irreducible 
$U_q(\widehat{ \mathfrak{sl}}_2)$-module
of type $(1,1)$.
We let $U_0, \ldots, U_d$
denote the 
weight space decomposition of $V$, 
from Lemma
\ref{lem:nnweight}.
We let $A:V\to V$
and $A^*:V\to V$ denote the linear
transformations that act as
$e_0^+ +K_0$
and $e_1^++K_1$ respectively.
By Lemma
\ref{lem:uqqs}, 
$V$ is an
${\mathcal A}_q$-module
on which
$x, y$ act as $A,A^*$ respectively.
\end{assumption}

\begin{remark}
\rm
\label{rem:caution}
Referring to Assumption
\ref{as}, the
${\mathcal A}_q$-module
 $V$ is not necessarily irreducible.
We will address this issue in 
Proposition
\ref{thm:imp}.
\end{remark}

\begin{lemma} 
\label{lem:acton}
With reference to
Assumption \ref{as}
the following hold 
for $0 \leq i \leq d$:
\begin{enumerate}
\item[{\rm (i)}]  
The element
$e^+_0$
acts on
$U_i$ as $A-q^{2i-d}I$.
\item[{\rm (ii)}]  
The element
$e^+_1$ 
acts on
$U_i$ as 
$A^*-q^{d-2i}I$.
\end{enumerate}
\end{lemma}
\noindent {\it Proof:} 
(i) 
The element $e^+_0$ acts on $U_i$ as $A-K_0$.
By
(\ref{eq:nnkmove}) and since 
$V$ has type $(1,1)$ we find
$K_0$ acts on $U_i$ as 
$q^{2i-d}I$. The result follows.
\\
\medskip
\noindent (ii) Similar to the proof of (i) above.
\hfill $\Box $ \\

\begin{lemma} 
\label{lem:acton2}
\label{eq:ar}
\label{eq:asl}
With reference to
Assumption \ref{as}
the following hold 
for $0 \leq i \leq d$:
\begin{enumerate}
\item[{\rm (i)}]  
$(A-q^{2i-d}I)U_i \subseteq U_{i+1}$,
\item[{\rm (ii)}]  
$(A^*-q^{d-2i}I)U_i \subseteq U_{i-1}$.
\end{enumerate}
\end{lemma}
\noindent {\it Proof:} 
(i) 
Combine
Lemma
\ref{lem:acton}(i) with the inclusion on
the left in
(\ref{eq:altemove1}).
\\
\noindent (ii)
Combine Lemma
\ref{lem:acton}(ii) with the inclusion on
the right in
(\ref{eq:altemove2}).
\hfill $\Box $ \\

\begin{lemma}
\label{lem:aeig}
With reference to 
Assumption \ref{as},
 each of $A,A^*$ is semisimple
with eigenvalues $q^{-d}, q^{2-d}, \ldots, q^d$.
Moreover for $0 \leq i \leq d$
the dimension of the eigenspace for $A$ (resp. $A^*$)
associated with $q^{2i-d}$ (resp. $q^{d-2i}$)
is equal to the dimension
of $U_i$.
\end{lemma}
\noindent {\it Proof:} 
We first display the eigenvalues of $A$.
Observe that the scalars 
$q^{2i-d}$ $(0 \leq i \leq d)$ are mutually
distinct since $q$ is not a root of unity.
Recall that $U_0, \ldots, U_d$
is a decomposition of $V$.
By Lemma \ref{lem:acton2}(i)
we see that, 
with respect to an appropriate basis for
$V$, $A$ is represented by a lower triangular
matrix that has diagonal entries
$q^{-d}, q^{2-d}, \ldots, q^d$, with $q^{2i-d}$ appearing
$\mbox{dim}(U_i)$ times for $0 \leq i \leq d$.
Hence for $0 \leq i \leq d$ the scalar
$q^{2i-d}$ is a root of the
characteristic polynomial of $A$ with multiplicity
$\mbox{dim}(U_i)$.
We now show $A$ is semisimple.
To do this we show that the minimal polynomial of
$A$ has distinct roots.
By Lemma \ref{lem:acton2}(i)
we find $\prod_{i=0}^d(A-q^{2i-d}I)$ vanishes on $V$.
By this and since
$q^{2i-d}$ $(0 \leq i \leq d)$ are distinct we see that
the minimal polynomial of $A$ has distinct roots.
Therefore $A$ is semisimple.
We have now proved our assertions concerning $A$; 
our assertions concerning $A^*$ are similarly proved.
\hfill $\Box $ \\

\begin{corollary}
\label{cor:nnil}
With reference to
 Assumption
\ref{as},
the ${\mathcal A}_q$-module
$V$
is NonNil.
\end{corollary}
\noindent {\it Proof:} 
Immediate from Lemma
\ref{lem:aeig}.
\hfill $\Box $ \\

\begin{definition}
\label{def:vi}
\rm
With reference to Assumption
\ref{as},
for $0 \leq i \leq d$ we let 
$V_i$ 
(resp. $V^*_i$)
denote the eigenspace of $A$
(resp. $A^*$)
associated
with the eigenvalue
$q^{2i-d}$ (resp. $q^{d-2i}$).
We observe that
$V_0,  \ldots, V_d$
(resp.
$V^*_0, \ldots, V^*_d$)
is a decomposition of $V$.
\end{definition}

\begin{lemma}
\label{lem:dsum}
With reference to 
Assumption \ref{as} 
and Definition
\ref{def:vi}
the following hold for
$0 \leq i \leq d$:
\begin{enumerate}
\item[{\rm (i)}]  
$V_i+ \cdots + V_d = U_i+ \cdots+ U_d$, 
\item[{\rm (ii)}]  
$V^*_0+ \cdots + V^*_i  = U_{0}+ \cdots+ U_i$. 
\end{enumerate}
\end{lemma}
\noindent {\it Proof:} 
\noindent (i) 
Let $X_i=\sum_{j=i}^d U_j$
and 
$X'_i=\sum_{j=i}^d V_j$.
We show $X_i=X'_i$.
Define $T_i = \prod_{j=i}^d (A-q^{2j-d}I)$.
Then $X'_i = \lbrace v \in V\,|\,T_iv=0\rbrace$,
and $T_iX_i=0$ by
Lemma \ref{lem:acton2}(i),
so $X_i \subseteq X'_i$.
Now define $S_i= \prod_{j=0}^{i-1}(A-q^{2j-d}I)$.
Observe that $S_iV=X'_i$, and
 $S_iV\subseteq X_i$ by
Lemma \ref{lem:acton2}(i),
so $X'_i \subseteq X_i$.
By these comments $X_i=X'_i$. 
\\
\noindent (ii) Similar to the proof of (i) above.
\hfill $\Box $ \\

\begin{lemma}
\label{lem:quasitd}
With reference to 
Assumption \ref{as}
and 
Definition
\ref{def:vi}
the following hold:
\begin{enumerate} 
\item[{\rm (i)}]  
 For $0 \leq i \leq d$,
\begin{eqnarray*}
A^*V_i \subseteq V_{i-1}+V_i + V_{i+1},
\end{eqnarray*}
where $V_{-1}=0$, $V_{d+1}=0$.
\item[{\rm (ii)}]  
For $0 \leq i \leq d$,
\begin{eqnarray*}
AV^*_i \subseteq V^*_{i-1}+V^*_i + V^*_{i+1},
\end{eqnarray*}
where $V^*_{-1}=0$, $V^*_{d+1}=0$.
\end{enumerate}
\end{lemma}
\noindent {\it Proof:} 
The elements $A,A^*$ satisfy the cubic
$q$-Serre relations by
Lemma \ref{lem:uqqs} and Assumption \ref{as}. The result
follows from this,
Lemma
\ref{lem:key}, and 
Definition 
\ref{def:vi}.
\hfill $\Box $ \\

\begin{remark}
\rm
With reference to Assumption \ref{as},
the pair $A,A^*$ is not necessarily a 
tridiagonal pair on $V$, since
$V$ might be reducible
as an 
${\mathcal A}_q$-module.
\end{remark}

\section{The projections $F_i, E_i, E^*_i$}

\begin{definition}
\label{def:f}
\label{def:e}
\label{def:es}
\rm
With reference to 
Assumption 
\ref{as} and Definition
\ref{def:vi}
we define the following for $0 \leq i \leq d$:
\begin{enumerate}
\item[{\rm (i)}]  
We let $F_i :V\to V$ denote the linear transformation
that satisfies
both
\begin{eqnarray*}
&& \qquad (F_i- I)U_i = 0, 
\\
&&F_i U_j = 0 \quad \mbox{if} \quad j \not=i, \qquad \qquad (0 \leq j \leq d).
\end{eqnarray*}
We observe $F_i$ is the projection from $V$ onto $U_i$.
\item[{\rm (ii)}]  
We let $E_i :V\to V$ denote the linear transformation
that satisfies
both
\begin{eqnarray*}
&& \qquad (E_i- I)V_i = 0, 
\\
&&E_i V_j = 0 \quad \mbox{if} \quad j \not=i, \qquad \qquad (0 \leq j \leq d).
\end{eqnarray*}
We observe $E_i$ is the projection from $V$ onto $V_i$.
\item[{\rm (iii)}]  
We let $E^*_i :V\to V$ denote the linear transformation
that satisfies
both
\begin{eqnarray*}
&& \qquad (E^*_i- I)V^*_i = 0, 
\\
&&E^*_i V^*_j = 0 \quad \mbox{if} \quad j \not=i, \qquad \qquad (0 \leq j \leq d).
\end{eqnarray*}
We observe $E^*_i$ is the projection from $V$ onto $V^*_i$.
\end{enumerate}
\end{definition}

\begin{proposition}
\label{lem:bij1}
With reference to
Assumption 
\ref{as} and
Definition \ref{def:e} the following hold
for $0 \leq i \leq d$:
\begin{enumerate}
\item[{\rm (i)}]  
The linear transformations
\begin{eqnarray*}
{{U_i\quad \rightarrow \quad  V_i}\atop {u \quad \rightarrow \quad E_iu}}
\qquad \qquad \qquad
{{V_i\quad \rightarrow \quad U_i}\atop {v \quad \; \rightarrow \quad F_iv}}
\end{eqnarray*}
are bijections, and  moreover, they are inverses.
\item[{\rm (ii)}]  
The linear transformations
\begin{eqnarray*}
{{U_i\quad \rightarrow \quad  V^*_i}\atop {u \quad \rightarrow \quad E^*_iu}}
\qquad \qquad \qquad
{{V^*_i\quad \rightarrow \quad U_i}\atop {v \quad \; \rightarrow \quad F_iv}}
\end{eqnarray*}
are bijections, and  moreover, they are inverses.
\end{enumerate}
\end{proposition}
\noindent {\it Proof:}
(i)
It suffices to show
$F_iE_i-I$ vanishes on $U_i$ and
$E_iF_i-I$ vanishes on $V_i$.
We will use the following notation.
For $0 \leq j \leq d$
recall $U_j+\cdots + U_d=V_j+\cdots + V_d$;
let $W_j$
denote this common sum.
We set $W_{d+1}=0$.
By the construction
 $W_i=U_i+W_{i+1}$ (direct sum)
and
$W_i=V_i+W_{i+1}$ (direct sum).
Also
$(I-F_i)W_i=W_{i+1}$
and 
$(I-E_i)W_i=W_{i+1}$.
We now show
$F_iE_i-I$ vanishes on $U_i$.
Pick $u \in U_i$.
Using $F_iE_i-I=(F_i-I)E_i+E_i-I$
and our preliminary comments
we routinely find
$(F_iE_i-I)u \in W_{i+1}$.
But 
 $(F_iE_i-I)u \in U_i$ by construction
 and $U_i\cap W_{i+1}=0$
 so
 $(F_iE_i-I)u=0$.
We now show
$E_iF_i-I$ vanishes on $V_i$.
Pick $v \in V_i$.
Using $E_iF_i-I=(E_i-I)F_i+F_i-I$
and our preliminary comments
we routinely find
$(E_iF_i-I)v \in W_{i+1}$.
But 
 $(E_iF_i-I)v \in V_i$ by construction
 and $V_i\cap W_{i+1}=0$
 so
 $(E_iF_i-I)v=0$.
We have now shown
$F_iE_i-I$ vanishes on $U_i$ and
$E_iF_i-I$ vanishes on $V_i$. Consequently
the given maps
are inverses. Each of these maps has an inverse
and is therefore a 
bijection.
\\
\noindent (ii) Similar to the proof of (i) above.
\hfill $\Box $ \\

\noindent The following formulae will be useful.

\begin{lemma}
\label{lem:e3}
With reference to
Assumption 
\ref{as} and
Definition \ref{def:e}, for $0 \leq i \leq d$ we have
\begin{eqnarray}
E_i &=& \prod_{{0 \leq j \leq d}\atop {j\not=i}}
\frac{A-q^{2j-d}I}{q^{2i-d}-q^{2j-d}},
\label{eq:eform}
\\
E^*_i &=& \prod_{{0 \leq j \leq d}\atop {j\not=i}}
\frac{A^*-q^{d-2j}I}{q^{d-2i}-q^{d-2j}}.
\label{eq:eforms}
\end{eqnarray}
\end{lemma}
\noindent {\it Proof:}
Concerning (\ref{eq:eform}),
let $E'_i$ denote the expression on the right in 
that line.
Using
Definition \ref{def:vi} we find
$(E'_i-I)V_i=0$
and $E'_iV_j=0$ $(0 \leq j \leq d, \; j\not=i)$.
By this and 
Definition
 \ref{def:e}(ii)
we find $E_i=E'_i$.
We have now proved
(\ref{eq:eform}). The proof of
(\ref{eq:eforms}) is similar.
\hfill $\Box $ \\

\section{How $E_0, E^*_0, P_V$ are related}

\noindent Our goal in this section is to prove the following
result.

\begin{theorem}
\label{lem:eep}
With reference to Assumption 
\ref{as} and
Definition
\ref{def:f},
 for 
$u \in U_0$ we have
\begin{eqnarray}
E^*_0E_0u = P_V(
q^{-1}(q-q^{-1})^{-2})u,
\label{eq:eep}
\end{eqnarray}
where  $P_V$ is from Definition
\ref{def:pv}.
\end{theorem}

\noindent We will use the following lemma.

\begin{lemma}
\label{lem:esumform}
With reference to Assumption 
\ref{as}
and Definition \ref{def:e} the following hold
for $0 \leq i,j\leq d$:
\begin{enumerate}
\item[{\rm (i)}]  
The action
of $E_i$ on $U_j$ is zero if $i<j$ and coincides with
\begin{eqnarray*}
\frac{1}
{(q^{2i-d}-q^{2j-d})(q^{2i-d}-q^{2j+2-d})
\cdots (q^{2i-d}-q^{2i-2-d})}
\end{eqnarray*}
times
\begin{eqnarray*}
\sum_{h=0}^{d-i}
\frac{(e^+_0)^{i-j+h}}
{(q^{2i-d}-q^{2i+2-d})
(q^{2i-d}-q^{2i+4-d})
\cdots
(q^{2i-d}-q^{2i+2h-d})}
\end{eqnarray*}
if $i\geq j$.
\item[{\rm (ii)}]  
The action
of $E^*_i$ on $U_j$ is zero if $i>j$ and coincides with
\begin{eqnarray*}
\frac{1}
{(q^{d-2i}-q^{d-2j})(q^{d-2i}-q^{d-2j+2})
\cdots (q^{d-2i}-q^{d-2i-2})}
\end{eqnarray*}
times
\begin{eqnarray*}
\sum_{k=0}^{i}
\frac{(e^+_1)^{j-i+k}}
{(q^{d-2i}-q^{d-2i+2})
(q^{d-2i}-q^{d-2i+4})
\cdots
(q^{d-2i}-q^{d-2i+2k})}
\end{eqnarray*}
if $i\leq j$.
\end{enumerate}
\end{lemma}
\noindent {\it Proof:} 
(i) Pick $u \in U_j$. We find $E_iu$.
First assume $i<j$. Observe $U_j\subseteq V_j+\cdots+V_d$
by Lemma
\ref{lem:dsum}(i) and $E_i$ is zero
on 
$V_j+\cdots+V_d$ so $E_iu=0$.
Next assume $i=j$.
By Lemma
\ref{lem:dsum}(i) 
and since $E_iu\in V_i$ we find
$E_iu\in U_i+\cdots+U_d$.
Consequently there exist
$u_s \in U_{s}$ $(i \leq s\leq d)$
such that $E_iu=\sum_{s=i}^{d} u_s$.
Using $(A-q^{2i-d}I)E_i=0$
we find
\begin{eqnarray}
0 &=& (A-q^{2i-d}I)E_iu
\nonumber
\\
&=&
(A-q^{2i-d}I)\sum_{s=i}^{d}{u_s}
\nonumber
\\
&=& \sum_{s=i}^{d}(e^+_0+q^{2s-d}-q^{2i-d})u_s
\label{eq:chaos}
\end{eqnarray}
in view of Lemma
\ref{lem:acton}(i).
Rearranging terms in
(\ref{eq:chaos}) we find
$0 =\sum_{s=i+1}^{d} u'_s$ where 
\begin{eqnarray*}
u'_s= e^+_0u_{s-1}+(q^{2s-d}-q^{2i-d})u_s
\qquad \qquad (i+1 \leq s \leq d).
\end{eqnarray*}
Since $u'_s \in U_s$ for
$i+1 \leq s \leq d$ and since
$U_0, \ldots, U_d$ is a decomposition we find
$u'_s=0$ for $i+1 \leq s \leq d$. Consequently
\begin{eqnarray*}
u_s=
(q^{2i-d}-q^{2s-d})^{-1}
e^+_0u_{s-1}
\qquad \qquad
(i+1\leq s\leq d).
\end{eqnarray*}
By Proposition
\ref{lem:bij1}(i) and since
$u_i=F_iE_iu$
we find
$u_i=u$. From these comments we get the result for $i=j$.
Now assume $i>j$.
Define
\begin{eqnarray}
v&=&\frac{(A-q^{2i-2-d}I)(A-q^{2i-4-d}I)\cdots (A-q^{2j-d}I)u}
{(q^{2i-d}-q^{2i-2-d})(q^{2i-d}-q^{2i-4-d})\cdots (q^{2i-d}-q^{2j-d})}.
\label{eq:uv}
\end{eqnarray}
Using $E_iA=q^{2i-d}E_i$
we find 
$E_iv=E_iu$. In order to get $E_iu$ we compute $E_iv$.
By Lemma
\ref{lem:acton}(i)
and Lemma \ref{lem:acton2}(i),
the numerator on the right in
(\ref{eq:uv}) is equal to $(e^+_0)^{i-j}u$
and contained in $U_i$.
Now 
 $v \in U_i$
in view of
(\ref{eq:uv}).
 We can now get 
$E_iv$ using our above discussion
concerning the case $i=j$.
By these comments we obtain the result for
$i>j$.
\\
\noindent (ii) Similar to the proof of (i) above.
\hfill $\Box $ \\

\noindent {\it Proof of Theorem 
\ref{lem:eep}:}
By Lemma
\ref{lem:esumform}(i) (with $i=0$, $j=0$) we find
\begin{eqnarray}
E_0u = \sum_{h=0}^d 
\frac{(e^+_0)^hu}
{(q^{-d}-q^{2-d})
(q^{-d}-q^{4-d})
\cdots
(q^{-d}-q^{2h-d})}.
\label{eq:e0sum}
\end{eqnarray}
Pick an integer $h$ $(0 \leq h \leq d)$.
By (\ref{eq:altemove1}) we have
$(e^+_0)^hu \in U_h$.
By this and Lemma
\ref{lem:esumform}(ii) (with $i=0$, $j=h$) we find
\begin{eqnarray}
E^*_0
(e^+_0)^hu =  
\frac{(e^+_1)^h(e^+_0)^hu}
{(q^{d}-q^{d-2})
(q^{d}-q^{d-4})
\cdots
(q^{d}-q^{d-2h})}.
\label{eq:e0ssum}
\end{eqnarray}
By Definition
\ref{def:sig},
\begin{eqnarray}
(e^+_1)^h(e^+_0)^hu=\sigma_h u.
\label{eq:sigmarole}
\end{eqnarray}
To obtain 
(\ref{eq:eep}) we multiply each term
in 
(\ref{eq:e0sum}) on the left
by $E^*_0$ and evaluate the results
using
(\ref{eq:e0ssum}),
(\ref{eq:sigmarole}), and
(\ref{eq:nbrack}),
(\ref{eq:nfac}),
(\ref{eq:pv}).
\hfill $\Box $ \\

\section{The shape of a tridiagonal pair}

\noindent 
In this section we return to the concept of a tridiagonal
pair, discussed in Section 2.
We will obtain a result concerning these objects that
will be of use in
our proof of Theorem
\ref{thm:2}. This result might be of independent interest.

\medskip
\noindent 
We will use the following notation.
Let $V$ denote a vector space over $\K$ with finite
positive dimension and let $A,A^*$ denote a tridiagonal
pair on $V$.
Let 
$V_0, \ldots, V_d$
(resp.
$V^*_0, \ldots, V^*_d$
)
denote a standard ordering of the eigenspaces of
$A$ 
(resp. $A^*$).
By
\cite[Corollary~5.7]{TD00}, for
$0 \leq i \leq d$ the spaces
$V_i, V^*_i$
have the same dimension;
we denote this common dimension by $\rho_i$.
By the construction $\rho_i \not=0$.
By
\cite[Corollary~5.7]{TD00}
and
\cite[Corollary~6.6]{TD00} the sequence
$\rho_0, \rho_1, \ldots, \rho_d$
is symmetric and unimodal; that is 
$\rho_i = \rho_{d-i}$ for $0 \leq i \leq d$
and 
$\rho_{i-1}\leq \rho_i$ for
$1 \leq i \leq d/2$.
We call the sequence 
$(\rho_0, \rho_1, \ldots, \rho_d)$
the {\it shape vector} of $A,A^*$.

\begin{theorem}
\label{thm:shape}
Let $V$ denote a vector space over $\K$ with finite positive
dimension and let $A,A^*$ denote 
a $q$-geometric
tridiagonal pair on $V$.
Let $(\rho_0, \rho_1, \ldots, \rho_d)$ denote the
corresponding 
shape vector.
 Then  there exists a nonnegative integer $N$ and
 positive integers $d_1, d_2, \ldots, d_N$ such that
\begin{eqnarray*}
\sum_{i=0}^d \rho_i z^i
 = \prod_{j=1}^N
(1+z+z^2+\cdots + z^{d_j}). 
\end{eqnarray*}
\end{theorem}
\noindent {\it Proof:} 
By Theorem
\ref{thm:maintd} 
there exists a
NonNil 
irreducible
${\mathcal A}_q$-module structure on $V$
such that the standard generators $x,y$ act as
$A,A^*$ respectively.
Let $(\alpha, \alpha^*)$ denote the
type of 
this ${\mathcal A}_q$-module.
Replacing 
$A,A^*$ by $A/\alpha, A^*/\alpha^*$ the type becomes
$(1,1)$ and the shape vector is unchanged.
By
Theorem 
\ref{thm:towardsuq}
there exists
an irreducible type $(1,1)$  
$U_q(\widehat{ \mathfrak{sl}}_2)$-module
structure on 
$V$ such that
$x$ and $y$ act as $e^+_0+K_0$  
and $e^+_1+K_1$ respectively.
Let $U_0, \ldots, U_d$ denote the corresponding
weight
space decomposition of $V$. 
By Lemma 
\ref{lem:aeig}
we find
$\rho_i = \mbox{dim}(U_i)$ for $0 \leq i \leq d$.
By this and
Lemma
\ref{cor:genf} we get the result.
\hfill $\Box $ \\

\section{A basis for 
${\mathcal A}_q$}

\noindent In this section we 
cite some results from
\cite{shape} that we will use in our proof of
Theorem
\ref{thm:2}.
We view 
the $\K$-algebra ${\mathcal A}_q$ as a vector space
over $\K$. We  display a basis for this vector space.

\medskip
\noindent 
Recall $x,y$ are the standard generators for 
${\mathcal A}_q$.

\begin{definition}
\rm
\label{def:redword}
Let $n$ denote a nonnegative integer. By a {\it word
of length $n$}
in 
${\mathcal A}_q$
 we mean
an expression of the form
\begin{eqnarray*}
a_1 a_2 \cdots a_n
\end{eqnarray*}
where $a_i = x$ or $a_i = y$ for $1 \leq i \leq n$.
We interpret the word of length 0 as the identity element
of 
${\mathcal A}_q$.
We say this word is {\it trivial}.
\end{definition}

\begin{definition}
\rm
Let $a_1 a_2 \cdots a_n$
denote a word in
${\mathcal A}_q$.
Observe there exists 
a unique sequence
$(i_1, i_2, \ldots, i_s)$
of positive integers such that
$a_1 a_2 \cdots a_n$ is one of
$
x^{i_1}y^{i_2}x^{i_3} \cdots y^{i_s} 
$
or
$
x^{i_1}y^{i_2}x^{i_3} \cdots x^{i_s} 
$
or 
$
y^{i_1}x^{i_2}y^{i_3} \cdots x^{i_s} 
$
or 
$
y^{i_1}x^{i_2}y^{i_3} \cdots y^{i_s} 
$.
We call the sequence
$(i_1, i_2, \ldots, i_s)$
the {\it signature}
of $a_1a_2 \cdots a_n$.
\end{definition}

\begin{example}
\rm
Each of the words $yx^2y^2x$, $xy^2x^2y$ has signature $(1,2,2,1)$.
\end{example}

\begin{definition}
\rm
Let $a_1 a_2 \cdots a_n$
denote a word in
${\mathcal A}_q$ and
let 
$(i_1, i_2, \ldots, i_s)$
denote the corresponding signature.
We say 
$a_1 a_2 \cdots a_n$
is {\it reducible}
whenever there exists an integer $\eta$ $(2 \leq \eta \leq s-1)$
such that $i_{\eta-1} \geq i_\eta < i_{\eta+1}$. We say a word
is {\it irreducible} whenever it is not reducible.
\end{definition}

\begin{example}
\rm
A word of length less than 4 is irreducible. The only reducible
words of length 4 are $xyx^2$ and $yxy^2$.
\end{example}

\noindent
In the following lemma we give a necessary and sufficient condition
for a given nontrivial word in 
${\mathcal A}_q$
to be irreducible.

\begin{lemma}
\label{lem:irredform}
{\rm \cite[Lemma~2.28]{shape}}
Let $a_1 a_2 \cdots a_n$
denote a nontrivial word in
${\mathcal A}_q$ and
let 
$(i_1, i_2, \ldots, i_s)$
denote the corresponding signature.
Then the following are equivalent:
\begin{enumerate}
\item[{\rm (i)}]  
The word 
$a_1 a_2 \cdots a_n$ is irreducible.
\item[{\rm (ii)}]  
There exists an integer $t$ $(1 \leq t \leq s)$ such
that
\begin{eqnarray*}
i_1 < i_2 < \cdots < i_t \geq i_{t+1} \geq i_{t+2}
\geq \cdots \geq i_{s-1} \geq i_s.
\end{eqnarray*}
\end{enumerate}
\end{lemma}

\begin{theorem}
\label{thm:irredw}
{\rm \cite[Theorem~2.29]{shape}}
The 
set of irreducible words in 
${\mathcal A}_q$ forms a basis for
${\mathcal A}_q$.
\end{theorem}

\section{From 
$U_q(\widehat{ \mathfrak{sl}}_2)$-modules to 
NonNil ${\mathcal A}_q$-modules, II }

\noindent In this section we return to the
situation of Assumption \ref{as}.
Referring to that assumption
let $W$ denote 
an irreducible ${\mathcal A}_q$-submodule of $V$.
Our next goal is to show that
$W$
contains
the space $V_0$
from Definition 
\ref{def:vi}. 



\begin{definition}
\label{def:ep}
\rm
With reference to Assumption \ref{as},
let $W$ denote an irreducible
${\mathcal A}_q$-submodule of $V$.
Observe that
$W$ is the direct sum of the nonzero
spaces among $E_0W,\ldots, E_dW$,
where the $E_i$ are from
Definition
\ref{def:e}(ii).
We define
\begin{eqnarray*}
r = \mbox{min} \lbrace i \,|\,0 \leq i \leq d, E_iW\not=0 \rbrace.
\end{eqnarray*}
We call $r$ the {\it endpoint} of  $W$.
\end{definition}

\begin{lemma}
\label{lem:dim1}
With reference to Assumption \ref{as},
let $W$ denote an irreducible
${\mathcal A}_q$-submodule of $V$ and let $r$ denote 
the endpoint of $W$.
Then $\mbox{dim}(E_rW)=1$.
\end{lemma}
\noindent {\it Proof:} 
By construction $W$
is an irreducible
${\mathcal A}_q$-module.
This module is NonNil
by
Lemma
\ref{lem:aeig}.
By this and Theorem
\ref{thm:maintd} the pair
$A|_W, A^*|_W$ is a $q$-geometric tridiagonal pair on $W$.
Let $s$ denote the
diameter 
of 
$A|_W,A^*|_W$.
Then $E_rW, E_{r+1}W, \ldots, E_{r+s}W$
is a standard ordering of the eigenspaces of
$A|_W$.
Applying
Theorem
\ref{thm:shape}
to
$A|_W, A^*|_W$ we find
$\mbox{dim}(E_rW)=1$.
\hfill $\Box $ \\

\begin{lemma}
\label{lem:luz}
With reference to 
Assumption \ref{as},
let $W$ denote an irreducible
${\mathcal A}_q$-submodule  of $V$
and let $r$ denote the endpoint of $W$.
Pick $v \in E_rW$ and  write $u=F_rv$.
Then $e^+_1 u = 0$.
\end{lemma}
\noindent {\it Proof:} 
Observe $u \in U_r$ by
Definition \ref{def:f}(i).
We assume $r\geq 1$; otherwise 
$e^+_1 u = 0$ since 
$e^+_1 U_0 = 0$.
Observe 
$e^+_1 u \in U_{r-1}$ by
(\ref{eq:altemove2}). In order to
show 
$e^+_1 u=0$ we show
$e^+_1 u \in U_r +\cdots + U_d$.
Since $A^*$ acts on $V$ as $e^+_1+K_1$
we have
\begin{eqnarray}
\label{eq:three}
e^+_1 u=A^*v-K_1 v + e^+_1(u-v).
\end{eqnarray}
We are going to show
that each of the three terms on the right
in 
(\ref{eq:three}) is contained in
$U_r +\cdots + U_d$.
By the definition of $r$ we have
$W = E_rW+ \cdots + E_dW$ 
so $W \subseteq V_r+\cdots + V_d$
in view of
Definition \ref{def:e}(ii).
By this and 
 Lemma
\ref{lem:dsum}(i) we find  
$W \subseteq U_r + \cdots + U_d$.
By construction $v \in W$ so
$A^*v \in W$. By these comments
$A^*v \in U_r +\cdots + U_d$.
We mentioned $v \in W$ so $v \in U_r + \cdots + U_d$.
Each of $U_r, \ldots, U_d$ is an eigenspace for
$K_1$ so $K_1v \in U_r + \cdots + U_d$.
Since 
$v
\in U_{r} + \cdots + U_d$
and since $u=F_rv$
we find
$u-v \in U_{r+1} + \cdots + U_d$.
Now $e^+_1 (u-v) \in U_r + \cdots + U_{d-1}$
by
(\ref{eq:altemove2}) so
$e^+_1 (u-v) \in U_r + \cdots + U_d$.
We have now shown that each of the three terms
on the right in
(\ref{eq:three}) is contained in
$U_r +\cdots + U_d$.
Therefore 
$e^+_1u \in U_r + \cdots + U_d$.
By this and since 
$e^+_1u \in U_{r-1}$ we find
$e^+_1u=0$.
\hfill $\Box $ \\

\begin{lemma}
\label{lem:updown}
With reference to Assumption \ref{as},
let $W$ denote an irreducible
${\mathcal A}_q$-submodule of $V$ and let $r$ denote 
the endpoint of $W$.
Pick $v \in E_rW$ and  write $u=F_rv$.
Then 
\begin{eqnarray}
(e^+_1)^i (e^+_0)^i u  \in \mbox{Span}(u)
\qquad \qquad (0 \leq i \leq d-r).
\label{eq:ud}
\end{eqnarray}
\end{lemma}
\noindent {\it Proof:}
We may assume $v \not=0$; otherwise the result is trivial.
Define
\begin{eqnarray}
\label{eq:bidef}
\Delta_i = 
(A^*-q^{d-2r}I)
(A^*-q^{d-2r-2}I)
\cdots
(A^*-q^{d-2r-2i+2}I).
\end{eqnarray}
Since $\Delta_i$ is a polynomial in
$A^*$ we find
$\Delta_iW\subseteq W$. In particular
$\Delta_iv \in W$ so
$E_r \Delta_iv \in E_rW$.
The vector $v$ spans
$E_rW$ by
Lemma
\ref{lem:dim1}
so
there exists $t_i \in \K$ such that
$E_r \Delta_iv = t_i v$.
By this and since $E_rv=v$ we find
$E_r (\Delta_i-t_iI)v=0$.
Now
$(\Delta_i-t_iI)v \in E_{r+1}W + \cdots + E_dW$
in view of
Definition \ref{def:ep}.
Observe
$E_{r+1}W + \cdots + E_dW \subseteq 
V_{r+1} + \cdots + V_d$
where the $V_j$ are from
Definition \ref{def:vi}.
By these comments and
Lemma 
\ref{lem:dsum}(i) we find  
$(\Delta_i-t_iI)v \in 
U_{r+1} + \cdots + U_d$.
Consequently
$F_r(\Delta_i-t_iI)v=0$. Recall $F_rv=u$ so
\begin{eqnarray}
F_r\Delta_i v=t_iu.
\label{eq:fbx}
\end{eqnarray}
We now evaluate $F_r\Delta_i v$.
Observe $v=E_ru$
by Proposition
\ref{lem:bij1}(i)
and since $u=F_r v$.
By Lemma
\ref{lem:esumform}(i)
(with $i=r$, $j=r$)
there exist nonzero scalars $\gamma_h\in \K$ $(0 \leq h\leq d-r)$
such
that
$v=\sum_{h=0}^{d-r}\gamma_h (e^+_0)^hu$.
For 
$0 \leq h \leq d-r$ 
we compute $F_r\Delta_i
(e^+_0)^hu$.
Keep in mind $(e^+_0)^hu \in U_{r+h}$
by 
the inclusion on the left in
(\ref{eq:altemove1}).
First assume $h<i$.
Using 
Lemma \ref{eq:asl}(ii) 
and 
(\ref{eq:bidef}) 
we find
$\Delta_i
(e^+_0)^hu$ 
is contained in $U_{r+h-i}+\cdots+ U_{r-1}$
so 
$F_r\Delta_i
(e^+_0)^hu=0$.
Next assume $h=i$.
By
Lemma
\ref{lem:acton}(ii)
and
(\ref{eq:bidef}) 
we find
$(\Delta_i-
(e^+_1)^i)
(e^+_0)^i u$ is contained in
$U_{r+1}+\cdots + U_{r+i}$.
By this and since
$(e^+_1)^i
(e^+_0)^i u \in U_r$ we find
$F_r\Delta_i
(e^+_0)^i u = 
(e^+_1)^i
(e^+_0)^i u$. 
Next assume $h>i$.
Using
Lemma \ref{eq:asl}(ii) 
and
(\ref{eq:bidef}) 
we find
$\Delta_i
(e^+_0)^hu$ 
is contained in $U_{r+h-i}+\cdots + U_{r+h}$.
By this and since $h>i$ we find
$F_r\Delta_i
(e^+_0)^h u =0$. 
By these comments
we find
$F_r\Delta_i v=
\gamma_i (e^+_1)^i
(e^+_0)^i u$.
Combining this and
(\ref{eq:fbx}) we obtain
(\ref{eq:ud}).
\hfill $\Box $ \\

\noindent We will use the following result.

\begin{lemma}
\label{lem:bt}
{\rm \cite[Theorem~13.1]{bt}}
Let $B$ denote the subalgebra of
$U_q(\widehat{ \mathfrak{sl}}_2)$ generated by
$e_0^+$, $e^+_1$, $K^{\pm 1}_0$, $K^{\pm 1}_1$.
Then with reference to Assumption \ref{as},
$V$ is irreducible as a $B$-module.
\end{lemma}

\begin{proposition}
\label{cor:eis0}
With reference to 
Assumption \ref{as},
let $W$ denote an irreducible
${\mathcal A}_q$-submodule  of $V$.
Then $W$ contains the space  $V_0$
from Definition \ref{def:vi}.
\end{proposition}
\noindent {\it Proof:}
We first
show $r=0$, where $r$ is the endpoint of $W$.
Pick a nonzero $v \in E_rW$ and write $u=F_rv$.
Observe 
$0 \not=u \in U_r$ by
Proposition
\ref{lem:bij1}(i).
By 
Lemma
\ref{lem:luz} and the inclusion on the left in
(\ref{eq:altemove1}),
\begin{eqnarray}
e_1^+u = 0,
 \qquad \qquad
(e_0^+)^{d-r+1}u=0.
\label{eq:ends}
\end{eqnarray}
By Lemma
\ref{lem:updown},
\begin{eqnarray}
(e^+_1)^i (e^+_0)^i u  \in \mbox{Span}(u)
\qquad \qquad (0 \leq i \leq d-r).
\label{eq:middle}
\end{eqnarray}
Let $\tilde W$ denote the subspace of
$V$ spanned by all vectors of the form
\begin{eqnarray}
(e^+_1)^{i_1}(e^+_0)^{i_2}(e^+_1)^{i_3}(e^+_0)^{i_4}
\cdots (e^+_0)^{i_n}u,
\label{eq:zigzag}
\end{eqnarray}
where $i_1, i_2, \ldots, i_n$ ranges over all sequences
such that $n$ is a nonnegative even integer,
and $i_1, i_2, \ldots, i_n$ are integers
satisfying
$0 \leq i_1 < i_2 <\cdots < i_n\leq d-r$.
Observe 
$u \in {\tilde W}$ so
${\tilde W}\not=0$.
We are going to show $\tilde W=V$
and  ${\tilde W}\subseteq U_r+\cdots+U_d$.
We now show
$\tilde W=V$.
To do this we show that
$\tilde W$ is 
$B$-invariant, where $B$ is from
Lemma
\ref{lem:bt}.
To begin, we show that 
$\tilde W$ is  invariant under
each of $e^+_0, e^+_1$.
By
Remark
\ref{rem:trad} there exists 
an 
${\mathcal A}_q$-module 
structure on $V$ such that
$x,y$ act as
$e^+_0, e^+_1$ respectively.
With respect to this
${\mathcal A}_q$-module structure
we have
$\tilde W= 
{\mathcal A}_q u$
in view of Lemma
 \ref{lem:irredform},
 Theorem
\ref{thm:irredw}
and 
(\ref{eq:ends}), 
(\ref{eq:middle}).
It follows that
$\tilde W$ is invariant under each of
$e^+_0, e^+_1$.
We now show  that 
$\tilde W$ is invariant under 
each of
$K^{\pm 1}_0, K^{\pm 1}_1$.
By 
(\ref{eq:altemove1}),
(\ref{eq:altemove2})
 the vector
(\ref{eq:zigzag}) is contained in $U_{r+i}$ where
$i=\sum_{h=1}^n i_h(-1)^h$.
Therefore the vector
(\ref{eq:zigzag}) 
is a common eigenvector for $K_0, K_1$.
Recall
$\tilde W$ is spanned by the vectors
(\ref{eq:zigzag}) 
so $\tilde W$ is invariant under 
each of
$K^{\pm 1}_0, K^{\pm 1}_1$.
By our above comments $\tilde W$ is
$B$-invariant. By this and
Lemma
\ref{lem:bt} we find $\tilde W=V$.
We now show
${\tilde W} \subseteq U_r + \cdots + U_d$.
We mentioned above that 
the vector (\ref{eq:zigzag}) is contained in $U_{r+i}$ where
$i=\sum_{h=1}^n i_h(-1)^h$.
From the construction $0\leq i \leq d-r$
so
$U_{r+i} \subseteq 
U_r+\cdots + U_d$.
Therefore 
the vector (\ref{eq:zigzag}) 
is contained in $U_r+\cdots + U_d$ so
${\tilde W}\subseteq U_r +\cdots + U_d$. 
We have shown
$\tilde W=V$ and
$\tilde W\subseteq 
 U_r + \cdots + U_d$.
Now 
$r=0$ since $U_0, \ldots, U_d$ is a decomposition
of $V$.
Now $E_0W\not=0$ by
Definition
\ref{def:ep}.
We have
$\mbox{dim}(V_0)=
\mbox{dim}(U_0)$
by Lemma
\ref{lem:aeig}
and 
$\mbox{dim}(U_0)=1$ by
Lemma
\ref{cor:genf}
so
$\mbox{dim}(V_0)=1$.
We have $0 \not= E_0W\subseteq V_0$ 
so 
$E_0W=V_0$.
But $E_0W\subseteq W$ by
(\ref{eq:eform})
so
$V_0\subseteq W$. 
%
\hfill $\Box $ \\


\section{The classification}

\noindent In this section we give 
the proof of
Theorem
\ref{thm:1} and
Theorem
\ref{thm:2}.
These proofs depend on the following
Proposition.

\begin{proposition}
\label{thm:imp}
With reference to Assumption \ref{as}
the following are equivalent:
\begin{enumerate}
\item[{\rm (i)}]  
$V$ is irreducible as an
${\mathcal A}_q$-module. 
\item[{\rm (ii)}]  
$P_V(q^{-1}(q-q^{-1})^{-2})\not=0$.
\end{enumerate}
\end{proposition}
\noindent {\it Proof:} 
(i) $\Rightarrow$ (ii)
We assume 
$P_V(q^{-1}(q-q^{-1})^{-2})=0$ and
get a contradiction.
Define
\begin{eqnarray*}
W_i = (V_0+\cdots + V_i)\cap (V^*_{i+1}+ \cdots + V^*_d)
\qquad \qquad (0 \leq i \leq d-1),
\end{eqnarray*}
where the $V_j, V^*_j$ are from Definition \ref{def:vi}.
Further define
$W=W_0+\cdots + W_{d-1}$.
We are going to show
that $W$ is
an
${\mathcal A}_q$-submodule of $V$  
and $W\not=V$,
$W\not=0$.
To begin, 
we first show
$AW\subseteq W$.
For 
 $0 \leq i \leq d-1$
we have $(A-q^{2i-d}I)\sum_{j=0}^i V_j
= 
\sum_{j=0}^{i-1} V_j$
by Definition
\ref{def:vi}
and $(A-q^{2i-d}I)\sum_{j=i+1}^d V^*_j
\subseteq 
\sum_{j=i}^d V^*_j
$
by Lemma 
\ref{lem:quasitd}(ii).
By these comments
\begin{eqnarray*}
&&(A - q^{2i-d}I)W_i \subseteq W_{i-1}
 \qquad (1 \leq i \leq d-1),
\qquad 
(A - q^{-d}I)W_0 = 0
\end{eqnarray*}
and it follows
$AW\subseteq W$.
We now show
$A^*W\subseteq W$.
For 
 $0 \leq i \leq d-1$
we have $(A^*-q^{d-2i-2}I)\sum_{j=0}^i V_j
\subseteq 
\sum_{j=0}^{i+1} V_j$
by Lemma 
\ref{lem:quasitd}(i)
and $(A^*-q^{d-2i-2}I)\sum_{j=i+1}^d V^*_j
=
\sum_{j=i+2}^d V^*_j
$
by Definition
\ref{def:vi}.
By these comments
\begin{eqnarray*}
(A^* - q^{d-2i-2}I)W_i \subseteq W_{i+1}
 \qquad (0 \leq i \leq d-2),
\quad 
(A^* - q^{-d}I)W_{d-1} = 0
\end{eqnarray*}
and it follows 
$A^*W\subseteq W$.
So far we have shown that 
$W$ is an
${\mathcal A}_q$-submodule  of $V$. 
We now show 
$W\not=V$.
For $0 \leq i \leq d-1$ we have
$W_i \subseteq V^*_{i+1} + \cdots + V^*_d$
so 
$W_i \subseteq V^*_1 + \cdots + V^*_d$.
It follows 
$W \subseteq V^*_1 + \cdots + V^*_d$ so
$W\not=V$.
We now show
$W\not=0$.
To do this we display a nonzero vector in $W_0$.
Pick a nonzero $u \in U_0$.
Applying Theorem
\ref{lem:eep} we find
$E^*_0E_0u=0$.
Write $v =E_0u$.
Then 
$0 \not=v 
\in 
 V_0$
by 
Proposition
\ref{lem:bij1}(i).
Observe $E^*_0v=0$ so
$v \in V^*_1 + \cdots + V^*_d$.
From these comments
$v \in W_0$.
We have displayed a nonzero vector $v$ contained in $W_0$.
Of course $W_0\subseteq W$  
so $W\not=0$.
We have now shown that 
$W$ is an
${\mathcal A}_q$-submodule of $V$
and
$W\not=V$,
$W\not=0$.
This contradicts our assumption that $V$ is 
irreducible as an
${\mathcal A}_q$-module.
We conclude
$P_V(q^{-1}(q-q^{-1})^{-2})\not=0$.
\\
\noindent 
(ii) $\Rightarrow$ (i)
Let $W$ denote an irreducible
${\mathcal A}_q$-submodule of $V$.
We show
$W=V$.
Define
$W'_i = 
W\cap U_i$ for
$0 \leq i \leq d$ and put
$W' = 
\sum_{i=0}^d W'_i$.
By the construction
$W'\subseteq W$.
We are going to show 
$W'=V$.
To do this we show
$W'\not=0$ and
$W'$ is $B$-invariant,
where $B$ is from
Lemma
\ref{lem:bt}.
We now show
$W'\not=0$.
To do this we display a nonzero vector in
$W'_0$.
By 
Proposition
\ref{cor:eis0} we find
$V_0 \subseteq W$.
Pick a nonzero $v \in V_0$
and write
$u=F_0v$. By
Proposition
\ref{lem:bij1}(i) we find
$0 \not=u \in U_0$
and $v=E_0u$.
Combining this last equation
with
Theorem
\ref{lem:eep} we find
$u= \alpha^{-1} E^*_0v$
where $\alpha=P_V(q^{-1}(q-q^{-1})^{-2})$.
We have $v \in W$ by construction 
and $E^*_0W \subseteq W$
by (\ref{eq:eforms}) so 
now $u \in W$.
By the above comments $0 \not= u \in W'_0$.
Of course $W'_0\subseteq W'$ so
$W'\not=0$.
We now show that $W'$ is $B$-invariant.
To begin, we show
$K_0^{\pm 1}W'\subseteq W'$.
By
(\ref{eq:nnkmove})
$K_0 -q^{2i-d}I$
vanishes on $U_i$ for $0 \leq i \leq d$.
Therefore
$K_0^{\pm 1} W'_i \subseteq W'_i$ for $0 \leq i \leq d$
so $K_0^{\pm 1}W'\subseteq W'$.
By this and
since 
$K_0K_1-I$ vanishes on $V$ we find
$K_1^{\pm 1}W'\subseteq W'$.
We show $e^+_0W' \subseteq W'$.
From 
Lemma
\ref{lem:acton}(i)
and the inclusion on the left in
(\ref{eq:altemove1}) we find
\begin{eqnarray*}
e^+_0W'_i \subseteq W'_{i+1}
\qquad  (0 \leq i \leq d-1),
\qquad 
e^+_0W'_d=0
\end{eqnarray*}
and it follows $e^+_0W'\subseteq W'$.
We show 
$e^+_1W' \subseteq W'$.
From 
Lemma
\ref{lem:acton}(ii)
and the inclusion on the right in
(\ref{eq:altemove2}) we find
\begin{eqnarray*}
e^+_1W'_i \subseteq W'_{i-1}
\qquad  (1 \leq i \leq d),
\qquad 
e^+_1W'_0=0
\end{eqnarray*}
and it follows $e^+_1W'\subseteq W'$.
By our above comments
$W'$ is $B$-invariant.
We have now shown that 
$W'$ is nonzero and $B$-invariant
so $W'=V$ by
Lemma
\ref{lem:bt}.
Recall $W'\subseteq W$ so $W=V$. The result follows.
\hfill $\Box $ \\

\medskip
\noindent
It is now a simple matter
to prove
Theorem
\ref{thm:1} and
Theorem
\ref{thm:2}.

\medskip
\noindent {\it Proof of Theorem 
\ref{thm:1}:}
Combine 
Theorem
\ref{thm:towardsuq}
and the implication
(i) $\Rightarrow$ (ii) in
Proposition
\ref{thm:imp}.
\hfill $\Box $ \\

\medskip
\noindent {\it Proof of Theorem 
\ref{thm:2}:}
By Corollary
\ref{cor:nnil} there
 exists a NonNil 
${\mathcal A}_q$-module 
structure on 
$V$ such that $x$ and  
$y$ act as $e^+_0+K_0$  
and $e^+_1+K_1$ respectively.
By construction this
${\mathcal A}_q$-module 
is
unique.
This
${\mathcal A}_q$-module 
is
irreducible by the implication
(ii) $\Rightarrow$ (i) in 
Proposition 
\ref{thm:imp}.
This
${\mathcal A}_q$-module 
is 
type $(1,1)$
by
Lemma \ref{lem:aeig}.
\hfill $\Box $ \\

\noindent Tatsuro Ito \hfil\break
\noindent Department of Computational Science \hfil\break
\noindent Faculty of Science \hfil\break
\noindent Kanazawa University \hfil\break
\noindent Kakuma-machi \hfil\break
\noindent Kanazawa 920-1192, Japan \hfil\break
\noindent email:  {\tt ito@kappa.s.kanazawa-u.ac.jp}

\bigskip

\noindent Paul Terwilliger \hfil\break
\noindent Department of Mathematics \hfil\break
\noindent University of Wisconsin \hfil\break
\noindent Van Vleck Hall \hfil\break
\noindent 480 Lincoln Drive \hfil\break
\noindent Madison, WI 53706-1388 USA \hfil\break
\noindent email: {\tt terwilli@math.wisc.edu }\hfil\break

\end{document}